\documentclass[a4paper]{amsart}

\usepackage{color}
\usepackage{graphicx,epstopdf}
\usepackage{amsmath}
\usepackage{amssymb}
\usepackage{amsfonts}
\usepackage{amsthm}

\theoremstyle{remark}

\begin{document}
 \title[]{Numerical study of fractional Camassa-Holm equations}

\author{Christian Klein$^{*}$}
\address{Institut de Math\'ematiques de Bourgogne, UMR 5584\\
Institut Universitaire de France\\
                Universit\'e de Bourgogne-Franche-Comt\'e, 9 avenue Alain Savary, 21078 Dijon
                Cedex, France\\
    E-mail Christian.Klein@u-bourgogne.fr}

\author{Goksu Oruc}
\address{Institut de Math\'ematiques de Bourgogne, UMR 5584\\
               Universit\'e de Bourgogne-Franche-Comt\'e, 9 avenue Alain Savary, 21078 Dijon
               Cedex, France\\
   E-mail topkarci@itu.edu.tr}
\date{\today}

\begin{abstract}
A numerical study of fractional Camassa-Holm equations is presented. 
Smooth solitary waves are constructed numerically. Their stability is 
studied as well as the long time behavior of solutions for general localised 
initial data from the Schwartz class of rapidly decreasing functions. 
The appearence of dispersive shock waves is explored. 
 
\end{abstract}

 
\thanks{This work was partially supported by 
 the ANR-17-EURE-0002 EIPHI, the Bourgogne
Franche-Comt\'e Region, the European fund FEDER,
 and by the 
European Union Horizon 2020 research and innovation program under the 
Marie Sklodowska-Curie RISE 2017 grant agreement no. 778010 
IPaDEGAN.  
}
\maketitle

\section{Introduction}
This paper is concerned with the numerical study of solutions to the  
fractional Camassa-Holm (fCH) equation given by
\begin{equation}\label{fCH}
u_t+ \kappa_1 u_x + 3 uu_x + D^{\alpha} u_t = -\kappa_2 [2D^{\alpha} 
(u u_x)+ uD^{\alpha}u_x ],
\end{equation}
where $\kappa_{1}$, $\kappa_{2}$ are real constants.
The fractional derivative $D^{\alpha}$ (also called fractional 
Laplacian) is defined via its Fourier 
symbol
\begin{equation}
	\mathcal{F}D^{\alpha} = |k|^{\alpha}
	\label{Da},
\end{equation}
where $\mathcal{F}$ denotes the Fourier transform and where $k$ is 
the dual Fourier variable, see (\ref{fourierdef}). 

\subsection{Background}

In the case $\alpha=2$, $\kappa_1= 2\omega$ and $\kappa_2=\frac{1}{3}$, the equation \eqref{fCH} turns into the Camassa Holm (CH) equation
\begin{equation}\label{ch}
u_t+ 2 \omega u_x  + 3 u u_x - u_{xxt} = 2 u_x u_{xx} + u u_{xxx}.
\end{equation}
The CH equation was first introduced in \cite{fokas} in a formal 
study of a class of integrable equations. In  \cite{camhol} the CH equation was 
presented to model unidirectional propagation of small-amplitude 
shallow water waves above a flat bottom. It has been also derived as 
geodesic flow on the circle in \cite{kolev02,kolev03,kouranbaeva} and 
recently in the context of nonlinear dispersive elastic waves in 
\cite{EEE}. The CH equation is completely integrable and has an infinite number of local conserved quantities \cite{fisher}, three of which are given in the following form:
\begin{equation*}
H_0=\int_{\mathbb{R}} udx, \hspace{15pt} H_1=\frac{1}{2} \int_{\mathbb{R}} (u^2 + u_x^2)dx,  \hspace{15pt}H_2=\frac{1}{2} \int_{\mathbb{R}} (u^3+ u u_x^2+ 2 \omega u^2)dx.
\end{equation*}
The CH equation has smooth solitary wave solutions and peaked 
solitons (\emph{peakons}) in the cases $\omega > 0$ and $\omega =0$, 
respectively. The existence of peaked solitary waves has been established in \cite{alber}.
A classification of travelling wave solutions that also contain 
cusped solitons (\emph{cuspons}) has been proposed in \cite{lenells}. 
The orbital stability  of the travelling waves has been obtained for 
smooth solitons in \cite{Constantin2002}, for peakons in 
\cite{Constantin2000} and for periodic peakons in 
\cite{lenells1,lenells2}. Wave breaking phenomena which cause 
solutions to remain bounded whereas their slopes blow up in finite 
time have been investigated in \cite{conswave}.  Many numerical 
approaches for the CH equation have been developed such as finite 
difference methods, finite-volume methods,  pseudo-spectral methods, 
local discontinuous Galerkin methods, see for instance 
\cite{finite1,finite2,finite3,artebrant,kalisch,galerkin,fengfeng,camassahuang,GK,AGK}.

Even though the CH equation has been studied extensively, there are 
only few results on the fractional CH equation which has been 
obtained in \cite{EEE} via a multiscales expansion for a fractional Boussinesq 
equation appearing in elasticity. In \cite{duruk} the local well-posedness of the Cauchy problem to the following form of fractional CH equation
\begin{equation}\label{fchv2}
u_t+ u_x + \frac{1}{2} (u^{2})_x + \frac{3}{4}D^{\alpha} u_x + \frac{5}{4}D^{\alpha} u_t = -\frac{1}{4}[2D^{\alpha} (u u_x)+ uD^{\alpha}u_x ]
\end{equation}
has been proven for the initial data $u_0\in$ $H^s(\mathbb{R})$, 
$s>\frac{5}{2}$ when $\alpha>2$ via Kato's semigroup approach for 
quasilinear evolution equations. In \cite{Besov} the local 
well-posedness criteria for the same Cauchy problem has been refined 
in an appropriate Besov space as $B_{2,1}^{s_0}$ with $s_0=2 \alpha - \frac{1}{2}$ for $\alpha > 3$ and $s_0=\frac{5}{2}$ for $2 < \alpha \leq 3$. In \cite{Besov2} the local well-posedness results given in \cite{Besov} have been recently extended to the Cauchy problem for the generalized fractional CH equation
\begin{equation}\label{gfchv2}
u_t+ u_x + \frac{1}{2} (u^{p+1})_x + \frac{3}{4}D^{\alpha} u_x + \frac{5}{4}D^{\alpha} u_t = -\frac{p+1}{8}[2D^{\alpha} (u^p u_x)+ u^pD^{\alpha}u_x ]
\end{equation}
with $\alpha>2$ and $p$ $\in$ $\mathbb{N}^+$.   A blow-up criterion 
for the solutions has been obtained. We note here that in the standard case ($\alpha=2$) the equation \eqref{fchv2} is related to well known integrable shallow water equation derived in \cite{dullin}. To the best of our knowledge there is no study related to fractional CH equation when $\alpha \leq 2$ and the current paper is the first numerical study on the fractional CH equation in the literature.

\subsection{Main results}
In this paper we always consider positive values of the constant 
$\kappa_{1}=2\omega$. For the CH equation it is known that the 
solitary waves with velocity $c>2\omega$ are always smooth, see 
\cite{Joh}. Our numerical study indicates that there might not be 
smooth solitary waves for positive $\omega$ and all values of 
$\alpha$. If $\alpha$ and thus the dispersion is too small, there 
might be no such solitary waves.  Thus it could be that for $\omega>0$, 
there exists a minimal value of $\alpha$ depending on $\omega$ and 
$c$, $\alpha_{s}(\omega,c)>0$, such that only for 
$\alpha>\alpha_{s}(\omega,c)$ there exist smooth solitary waves. It also 
appears that for given 
$\alpha$ and $\omega>0$, there exists a minimal velocity 
$c_{s}(\omega,\alpha)$ such 
that there exist smooth solitary waves for $c>c_{s}$. Since we 
construct the solitary waves numerically with a Newton iteration, the 
failure of the iteration only gives an indication that there are no 
smooth solutions for certain parameters $\alpha$, $c$, $\omega$. 
However, this could also mean that there are no appropriate initial 
iterates known. 

By studying perturbations of the numerically constructed solitary 
waves, we get \\
\textbf{Main conjecture I:}\\
The smooth solitary waves are 
orbitally stable.

For the CH equation,  it was shown in \cite{McK}, that solutions for 
smooth initial  
data $u_{0}$ subject to the condition 
$u_{0}-\partial_{xx}u_{0}+\omega>0$ will stay smooth for all times. 
The precise nature of the singularity that can appear in finite time 
for initial data not satifying this non-breaking condition does not 
appear to be known. Our numerical experiments indicate that for 
sufficiently small $\alpha$, initial data satisfying a  condition of 
the form $u_{0}+D^{\alpha}u_{0}+\omega>0$ can develop a cusp in 
finite time. We have\\
\textbf{Main conjecture II:}\\
For sufficiently small $\alpha<\alpha_{c}(\omega)$, initial data of 
sufficiently large mass lead to a blow-up of the fCH solution in 
finite time. The blow-up near a point $x_{s}$ is a cusp of the form 
$u\propto\sqrt{|x-x_{s}|}$.

The precise value of $\alpha_{c}(\omega)$ and the conditions on the 
initial data are not known.

This paper is organized as follows: in Section 2 basic facts on the 
standard CH and the fractional CH equation are gathered and some 
useful notation is reviewed. In Section 3 the numerical approach is 
introduced for solitary waves to the fractional CH equation. In 
Section 4 we investigate numerically the stability properties of 
solitary waves. The long time behavior of solution to fractional CH 
equation is studied in the case of initial data from the Schwartz class 
in Section 5. In Section 6, we study the appearence of rapid 
modulated oscillations, \emph{dispersive shock waves}, in the 
vicinity of shocks to the corresponding dispersionless equation. We 
add some concluding remarks in Section 7. 

\section{Preliminaries}
In this section we collect some basic facts on the standard CH and 
the fractional CH equation. 

We apply the standard definition for a  Fourier transform for  
tempered distributions $u(x)$
denoted by $(\mathcal{F} u)(k)$ with dual variable $k$ 
and its inverse,
\begin{align}
	(\mathcal{F}u)(k)=\hat{u}& = \int_{\mathbb{R}}^{}u(x)e^{-ikx}dx,\quad k\in 
	\mathbb{R},
	\nonumber\\
	u(x) & =\frac{1}{2\pi}\int_{\mathbb{R}}^{}(\mathcal{F}u)(k)e^{ikx} 
	dk,\quad x\in\mathbb{R}.
	\label{fourierdef}
\end{align}

\subsection{Conserved Quantities}
We give the derivation of the conserved quantities for the equation 
\eqref{fCH}. To this end we consider sufficiently smooth  solutions 
which tend to $0$ as $x\rightarrow \mp \infty$. Integrating equation 
\eqref{fCH} over the real line, we get
\begin{equation} \label{fCH1}
  \frac{d}{dt} \int_{\mathbb{{R}}} \left( I+  D^{\alpha}  \right) u dx +   \int_{\mathbb{{R}}}\left(\kappa_1 u + \frac{3}{2} u^2  + \kappa_2{D^{\alpha}u^2} \right)_x dx + \kappa_2 \int_{\mathbb{{R}}} u D^{\alpha}u_x dx =0.
\end{equation}
By using Plancherel's theorem, we can rewrite the last integral in equation \eqref{fCH1} as
\begin{eqnarray*}
  \int_{\mathbb{{R}}} u(x,t) D^{\alpha}u_x(x,t) dx  &=& 
  \int_{\mathbb{{R}}} \hat{u}(k,t) |k|^{\alpha} \overline{ i k \hat{u}(k,t)} \frac{dk}{2\pi}  , \nonumber \\
   &=& \int_{\mathbb{{R}}} |k|^{\frac{\alpha}{2}} \hat{u} 
   |k|^{\frac{\alpha}{2}}  \overline{ i k \hat{u}} \frac{d k}{2\pi}  , \nonumber \\
   &=&   \int_{\mathbb{{R}}}
   D^{\frac{\alpha}{2}} u D^{\frac{\alpha}{2}} u_x dx,\nonumber \\
   &=&   \frac{1}{2} \int_{\mathbb{{R}}}
   |D^{\frac{\alpha}{2}} u(x,t)|_x^2 dx.
\end{eqnarray*}
This implies for equation \eqref{fCH1}
\begin{equation} \label{fCH2}
  \frac{d}{dt} \int_{\mathbb{{R}}} \left( I+  D^{\alpha}  \right) u dx +   \int_{\mathbb{{R}}}\left(\kappa_1 u + \frac{3}{2}u^2 + \kappa_2 {D^{\alpha}u^2} +\frac{\kappa_2}{2} |D^{\frac{\alpha}{2}} u|^2 \right)_x dx =0,
\end{equation}
which gives the conserved mass of the fCH equation in the following form
\begin{equation}\label{mass}
  I_1= \int_{\mathbb{{R}}} \left( u(x,t)+ D^{\alpha} u(x,t) \right) dx.
\end{equation}

For the second conserved quantity, we first multiply both sides of 
the fCH equation by $u$ and integrate over  $\mathbb{{R}}$. Then we have
\begin{equation*}
\frac{1}{2} \frac{d}{dt} \int_{\mathbb{{R}}} \left( u^2+ |D^{\frac{\alpha}{2}}u|^2  \right)  dx +   \int_{\mathbb{{R}}}\left( \frac{\kappa_1}{2}u^2 + {u^3} +\kappa_2 u^2 D^{\alpha}u \right)_x dx=0.
\end{equation*}
Here we have used 
\begin{equation*}
\int_{\mathbb{{R}}} u(x,t) D^{\alpha}(u^2(x,t))_x dx  =  \int_{\mathbb{{R}}}
D^{\alpha} u(x,t) (u^2(x,t))_x  dx.
\end{equation*}
Finally, the following identity is obtained as formal conserved energy of fCH equation
\begin{equation}\label{energy}
	 I_2= \int_{\mathbb{{R}}} \left(u^2(x,t) + |D^{\frac{\alpha}{2}}u(x,t) |^2 \right) dx.
\end{equation}

\subsection{Solitary Waves of the CH equation}
Solitary waves are localised traveling waves, i.e., solutions of the 
form
$u(x,t)=Q_c(\xi), ~~\xi=x-ct$ with constant propagation speed $c$ and 
fall-off condition $\displaystyle{\lim_{|\xi| \rightarrow \infty} 
Q_c(\xi)=0}$. This ansatz leads to equation (\ref{sw0}) for the fCH 
equation. 

For the integrable CH equation ($\alpha=2$, $\kappa_{1}=2\omega$, 
$\kappa_{2}=1/3$), equation (\ref{sw0})) can be integrated explicitly leading to 
\begin{equation}
	(-c(1-\partial_{xx})+2\omega)Q_{c}+\frac{3}{2}Q_{c}^{2}=Q_{c}Q_{c}''+\frac{1}{2}(Q_{c}')^{2}.
	\label{QCH}
\end{equation}
As discussed in \cite{Joh}, the soliton is given implicitly by 
\begin{equation}
	Q_{c} = 
	\frac{(c-2\omega)\mbox{sech}^2(\theta)}{\mbox{sech}^2(\theta)+2(\omega/c)
		\tanh^2(\theta)}
	\label{QCH2}
\end{equation}
where 
\begin{equation}
	x-ct = 
	\sqrt{\frac{4c}{c-2\omega}}\theta+\ln\frac{\cosh(\theta-\theta_{0})}{\cosh(\theta+\theta_{0})}
	\label{QCH3}
\end{equation}
with $\theta_{0}=\mbox{arctanh}(\sqrt{1 - (2\omega/c)})$. There is a 
smooth soliton for $\omega>0$ and $c>2\omega$. These solitons are 
orbitally stable, see \cite{Joh}.

\section{Solitary Waves}
In this section we numerically construct  solitary waves of the 
fractional CH equation. 

\subsection{Defining equations}

With the traveling wave ansatz 
$u(x,t)=Q_c(\xi), ~~\xi=x-ct$, 
equation \eqref{fCH} reduces to the following equation
\begin{equation}\label{sw0}
  -cQ_c' + \kappa_1 Q_c' + \frac{3}{2}(Q_c^{2})'  -cD^{\alpha} Q_c' +  \kappa_2 [ D^{\alpha} (Q_c^2)' + Q_c D^{\alpha} Q_c' ]  =0.
\end{equation}
Here $^\prime$ denotes the derivative with respect to $\xi$.
Integrating we get with the fall-off condition at infinity
\begin{equation}
	(-c(1+D^{\alpha})+\kappa_{1})Q_{c}+\frac{3}{2}Q_{c}^{2}+\kappa_{2}(D^{\alpha}Q_{c}^{2}+
	\partial_{\xi}^{-1}(Q_{c}D^{\alpha}Q_{c}'))=0.
	\label{Q}
\end{equation}
Here the antiderivative is defined via its Fourier symbol,  
$\mathcal{F}\partial_{x}^{-1}=1/(ik)$, i.e.,
$\partial_{x}^{-1}=\frac{1}{2}(\int_{-\infty}^{x}-\int_{x}^{\infty})$.

\subsection{Numerical approach}
To construct solitary waves numerically, we apply as in \cite{KS15} a Fourier spectral 
method. We study equation (\ref{Q}) in the Fourier domain  
\begin{equation}
	(-c(1+|k|^{\alpha})+\kappa_{1})\mathcal{F}Q_{c}+\frac{3}{2}\mathcal{F}(Q_{c}^{2})+\kappa_{2}(|k|^{\alpha}\mathcal{F}(Q_{c}^{2})
	+\frac{1}{ik}\mathcal{F}(Q_{c}D^{\alpha}Q_{c}'))=0.
	\label{Qfourier}
\end{equation}
The Fourier transform is approximated on a sufficiently large torus 
($x\in L[-\pi,\pi]$ with $L\gg1$) via a discrete Fourier transform 
(DFT) which is conveniently 
computed by a fast Fourier transform (FFT). This means we introduce 
the standard discretisation $x_{n}=-\pi L+nh$, $n=1,\ldots,N$, 
$h=2\pi L/N$ with $N\in \mathbb{N}$ the number of Fourier modes. The 
dual variable is then given by $k=(-N/2+1,\ldots,N/2)/L$. In an abuse 
of notation, we will use the same symbols for the DFT as for the 
standard Fourier transform.

With the DFT discretisation, equation (\ref{Qfourier}) becomes a 
system of $N$ nonlinear equations for $\mathcal{F}Q_{c}$ which is 
solved with a Newton-Krylov method. This means that the action of the 
inverse of the Jacobian in a standard Newton method is computed 
iteratively with the Krylov subspace method GMRES \cite{SS86}. Since the 
convergence of a Newton method is local, the choice of the initial 
iterate is important. Therefore we apply a tracing technique: the 
implicit solution for CH ($\alpha=2$) for some value of the velocity 
$c$ is taken as the initial iterate for some smaller value of 
$\alpha$. The result of this iteration is taken as an initial iterate 
of an even smaller value of $\alpha$. 

Numerically challenging is the computation of the last term in 
(\ref{Qfourier}) since there is a division by the dual variable $k$ 
that can vanish. 
We apply the same approach as a in \cite{KMS}, the limiting value for $k=0$ is 
computed via de l'Hospital's rule, $\lim_{k\to0} \hat{f}(k)/k = \hat{f}'(0)$ if 
$\hat{f}'(k)$ is bounded for $k=0$. The derivative for $k=0$ is computed in 
standard way via the sum of the inverse Fourier transform $f$ of the 
function $\hat{f}$ times $x$ sampled at the 
collocation points in $x$, see \cite{KMS}, i.e., $\hat{f}'(0) \approx
\sum_{n=1}^{N}x_{n}f_{n}$. 

\subsection{Examples}

\begin{figure}[htb!]
 \includegraphics[width=\textwidth]{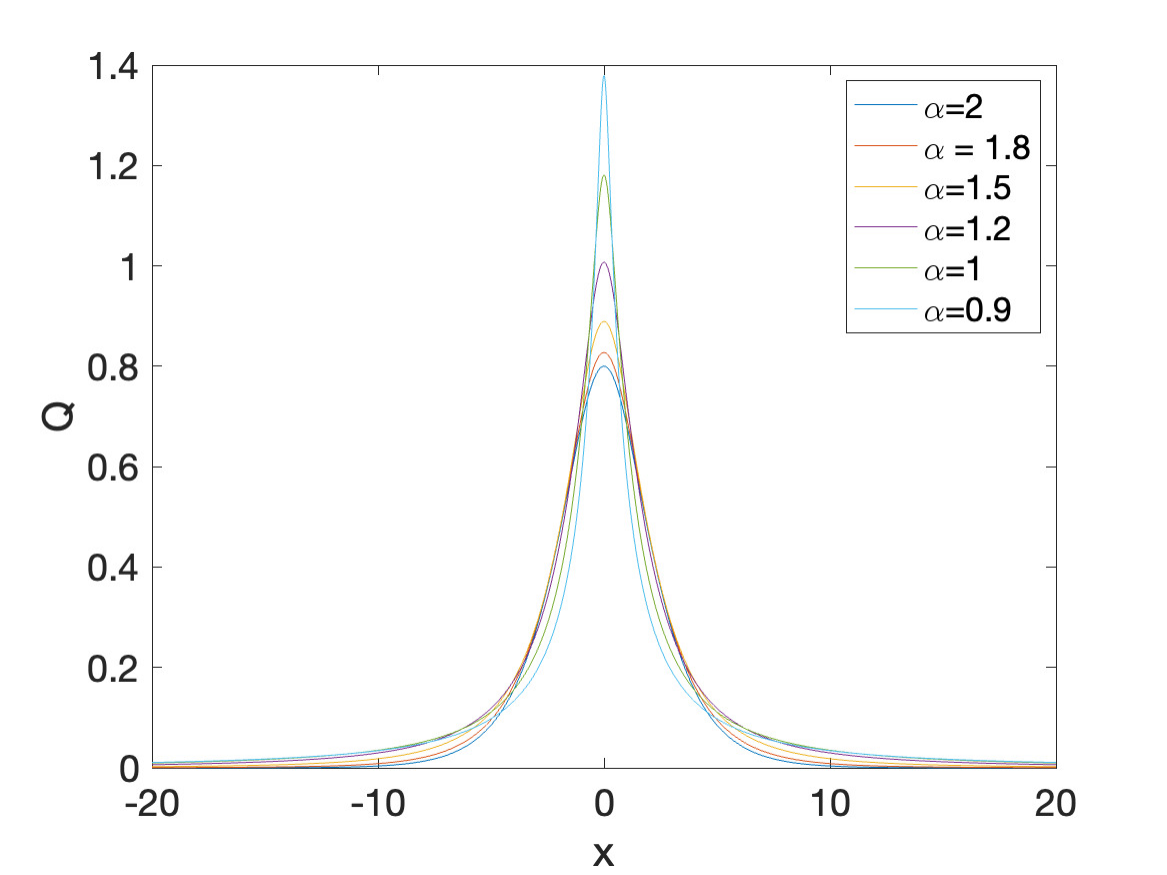}
 \caption{Solitary waves  (\ref{Q}) for $c=2$, $\kappa_{1}=1.2$, 
 $\kappa_{2}=1/3$ for several values of $\alpha$. }
 \label{Qc2}
\end{figure}

To study concrete examples, we work with $L=100$ and $N=2^{16}$ as in 
\cite{KS15} for fractional Korteweg-de Vries (fKdV) equations. We consider the case 
$\kappa_{2}=1/3$, $\kappa_{1}=2\omega$ being integrable for 
$\alpha=2$. We choose $\omega=3/5$ as in \cite{EEE}. We show solitary 
waves for $c=2$ for several values of $\alpha$ in Fig.~\ref{Qc2}. It 
can be seen that the smaller $\alpha$ and thus the dispersion, the 
more the solitary wave is peaked and the slower the decay towards 
infinity. We were not able to numerically find a solution for even 
smaller values of $\alpha$ since the iteration did not converge to a 
smooth solution. The behavior of the DFT coefficients in the 
iteration indicates that it converges to a peakon or cuspon. 
This does not mean that there are no smooth solitary 
waves for smaller values of $\alpha$ for this velocity, we just could 
not find them with a Newton iteration. But this could indicate that 
there might be for a given $c$ a lower limit 
$\alpha_{s}(\omega,c)$ for $\alpha$ below which there are no 
smooth solitons; we recall that there are smooth solitons for CH 
for all positive $\omega$ with $c>2\omega$.  
Note, however, that it is difficult to identify such 
a limit with an iterative method since the failure of the latter to 
converge does not mean that there is no such solution. It can be that 
the initial iterate was just not sufficiently close. 

\begin{figure}[htb!]
 \includegraphics[width=0.49\textwidth]{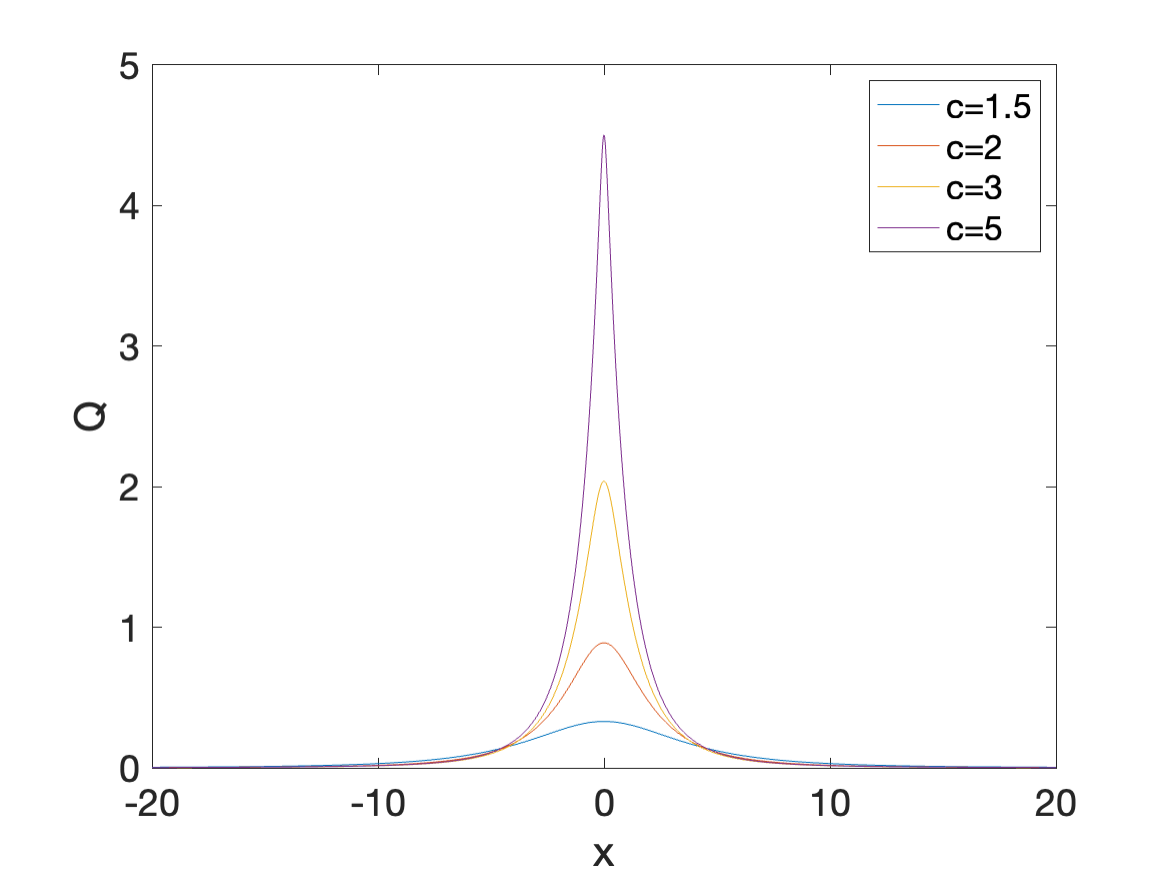}
 \includegraphics[width=0.49\textwidth]{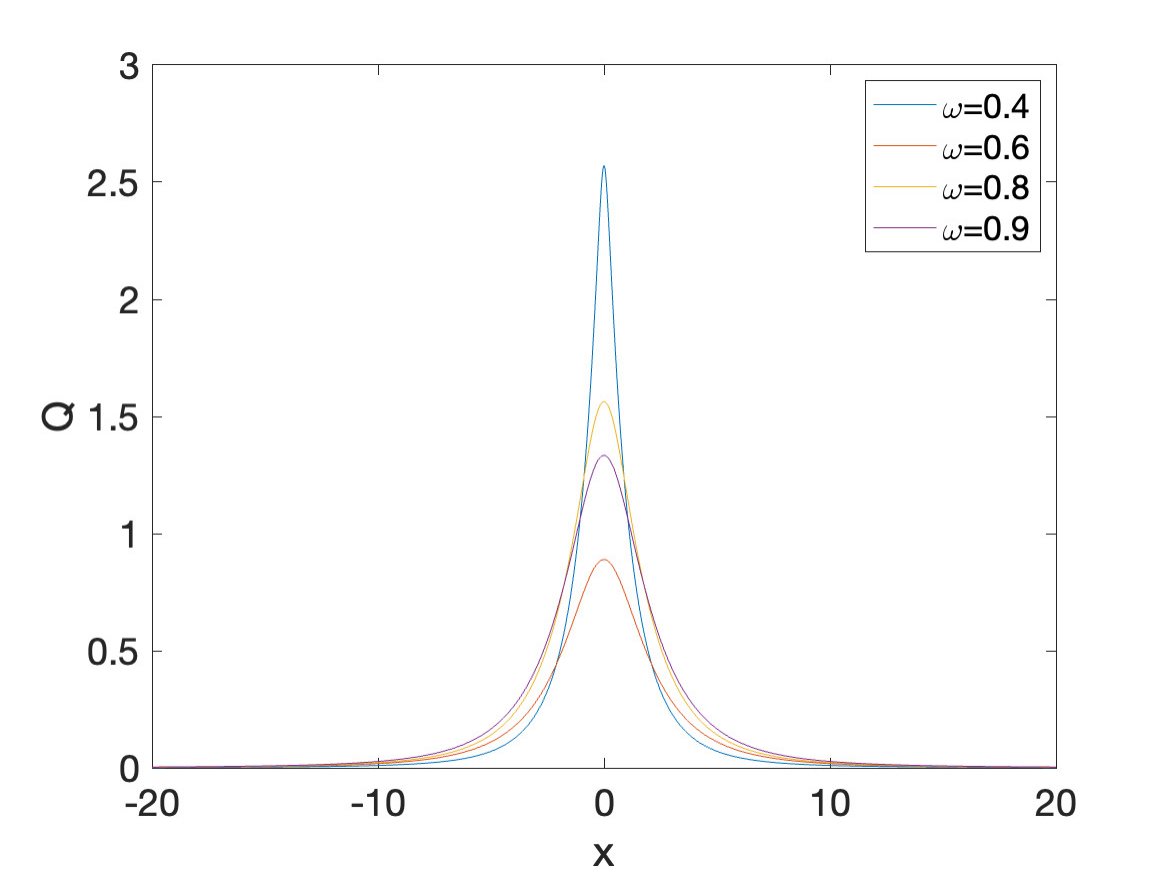}
 \caption{Solitary waves  (\ref{Q}) for  
 $\kappa_{2}=1/3$ and $\alpha=1.5$, on the left for $\omega=0.6$ and 
 several values of $c$, on the right for $c=2$ and several values of $\omega$. }
 \label{Qa15}
\end{figure}

The amplitude of the solitary waves increases with the velocity $c$. 
We show this behavior on the left of Fig.~\ref{Qa15} for $\alpha=1.5$, 
$\omega=0.6$ and several valus of $c$. For CH solitons one has 
$c>2\omega$. There is clearly also a lower bound on the velocity for 
smooth solitons larger than the $c>2\omega$ limit for CH, but it is 
not clear whether it is larger than the value for CH. Once more such 
a value cannot be determined with an 
iterative approach. If one fixes $c$ and $\alpha$,  there does not appear to be an obvious dependence 
between $\omega$ and the amplitude as in the case 
of CH as can be seen on the right of Fig.~\ref{Qa15}. 

Note that the precise fall-off of the solitary waves for 
$|x|\to\infty$ appears to be unknown, it should be the same $|x|^{-(1+\alpha)}$ 
rate of solitary waves of the 
fractional Kortweg-de Vries and nonlinear Schr\"odinger equation, see 
\cite{FL}. But it is numerically difficult to determine the exact 
fall-off rate on the real line via an approximation on a torus.  

\section{Stability of the solitary waves}
In this section we study the stability of the solitary waves under 
small perturbations. For various perturbations, the solitary waves 
appear to be stable providing numerical evidence to the first part of 
the Main Conjecture. 

\subsection{Numerical approach}
To study the time evolution of solutions to the fractional CH equation 
(\ref{fCH}), we use the same spatial discretisation as in the 
previous section, i.e., a standard FFT approach. The time integration 
is done with the well known explicit Runge-Kutta method of 4th order. 
The accuracy of the time integration is controlled via the conserved 
energy (\ref{energy}) which will numerically depend on time due to 
unavoidable numerical errors. As discussed in \cite{etna}, this 
numerically computed energy will overestimate the numerical error by 
2-3 orders of magnitude. Thus we will always track in the following 
the relative energy and the DFT coefficients to control the numerical 
resolution. 

As a test we propagate the solitary wave for $c=2$, $\alpha=1.5$, 
$\omega=3/5$ and use $N_{t}=10^{4}$ time steps for $t\in[0,1]$. The 
relative energy is conserved during the whole computation to the 
order of $10^{-15}$. The DFT coefficients decrease to machine 
precision (here $10^{-16}$) which means that the solution is well resolved both in 
space and in time. The difference between the solitary wave 
travelling with velocity $c=2$ and the numerically computed solution 
for $t=1$ is of the order $10^{-15}$. This shows that the code is 
able to propagate solitary waves with machine precision as indicated 
by the numerical conservation of the energy and the DFT coefficients. 
In addition it shows that the solitary waves numerically constructed 
in the previous section are indeed solutions to equation (\ref{Q}) 
with a similar accuracy. 

\subsection{Perturbed solitary waves}
We consider perturbations of the solitary waves in the form of 
perturbed initial data,
\begin{equation}
	u(x,0) = Q_{c}+A e^{-x^{2}},
	\label{gausspert}
\end{equation}
where $A$ is a small real constant. In all examples below, the DFT 
coefficients always decrease to machine precision, and the relative 
energy is conserved to better than $10^{-6}$. 

\begin{figure}[htb!]
 \includegraphics[width=0.49\textwidth]{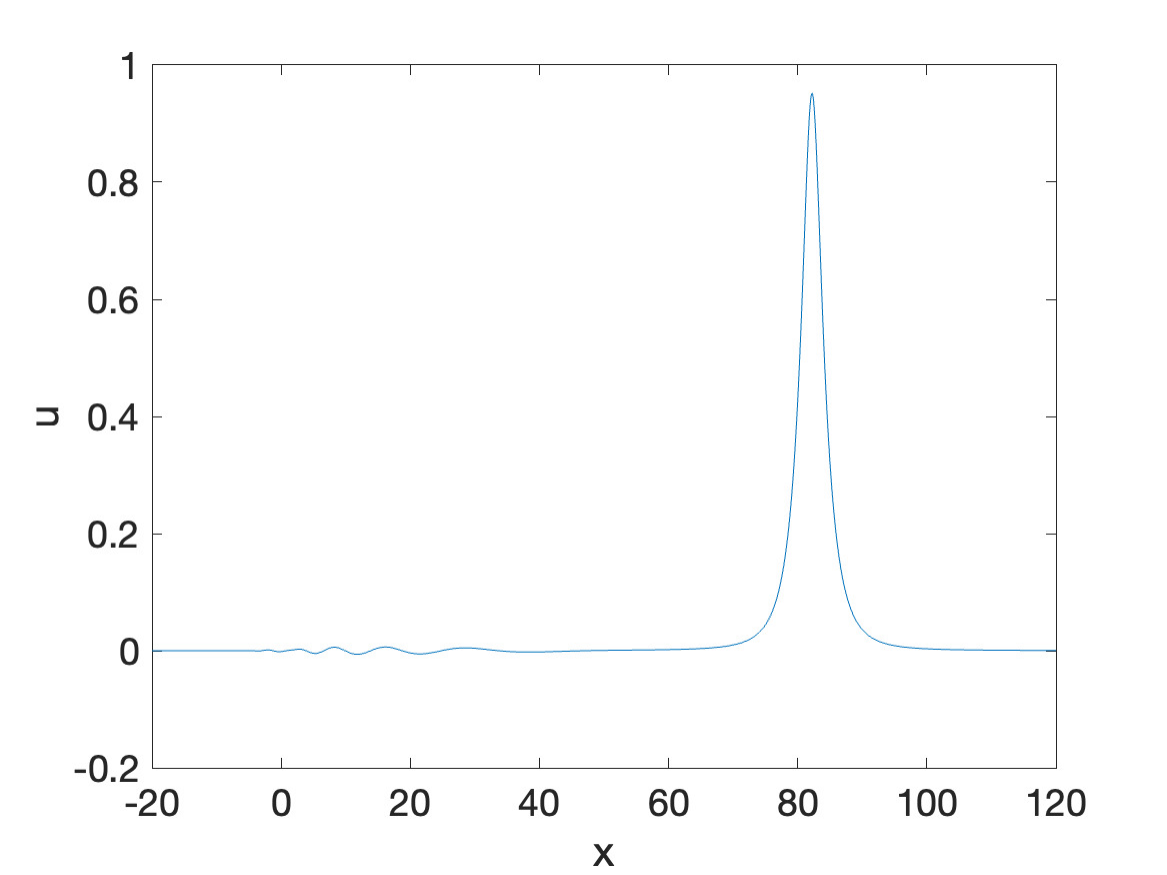}
 \includegraphics[width=0.49\textwidth]{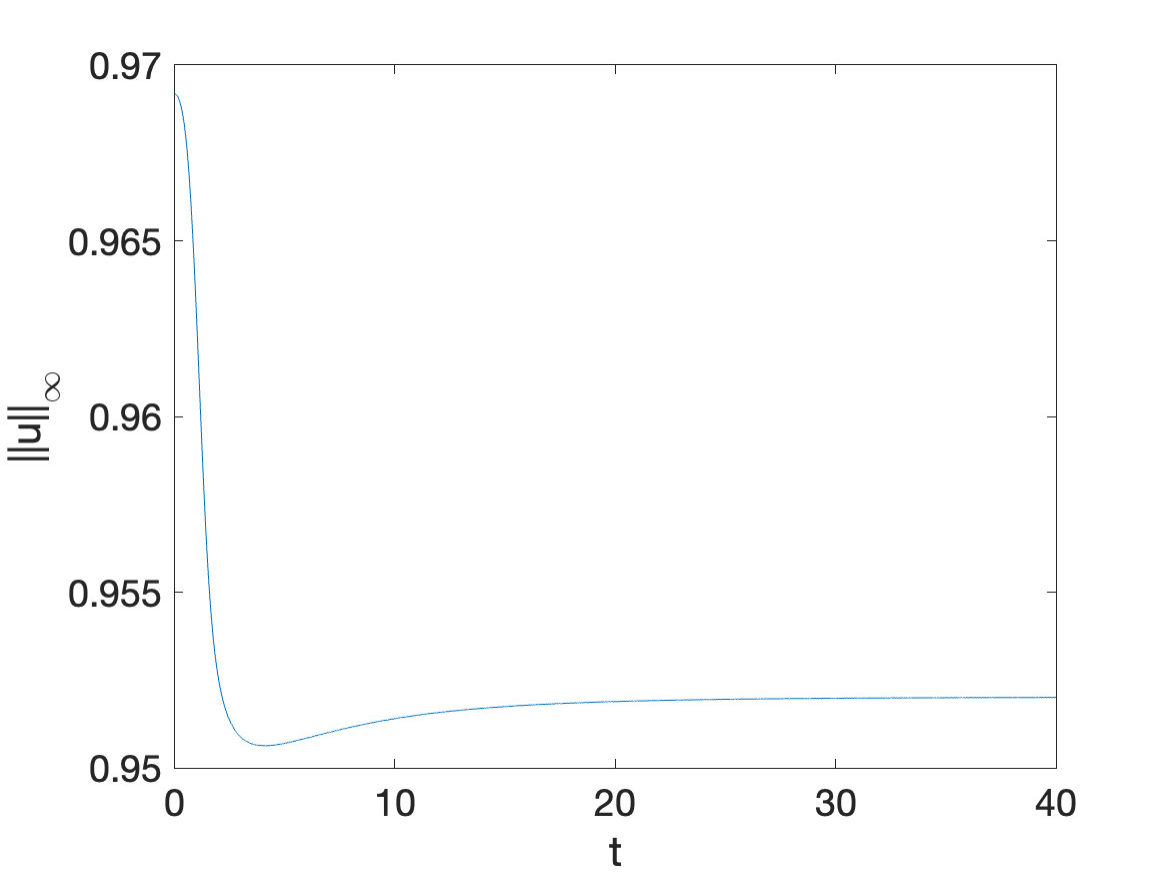}
 \caption{Solution to the fractional CH equation for initial data of 
 the form (\ref{gausspert}) with $A=0.08$ for $\alpha=1.5$ and 
 $c=2$: on the left the solution for $t=40$, on the right the 
 $L^{\infty}$ norm of the solution in dependence of time. }
 \label{Qa15gausspert}
\end{figure}

First we consider the case $\alpha=1.5$ for $Q_{2}$ and $A=\pm 0.08$. 
This corresponds to a perturbation of the order of 10\%. This is not 
a small perturbation, but in order to see numerical effects of a 
perturbation in finite time, it is convenient to consider 
perturbations that are of the order of a few percent. We use 
$N_{t}=10^{4}$ time steps for $t\leq40$. In Fig.~\ref{Qa15gausspert} we 
show on the left the solution for $t=40$. The solution appears to be 
a solitary wave with some radiation propagating towards $-\infty$. 
This is confirmed by the $L^{\infty}$ norm of the solution on the 
right of the same figure. After some time the $L^{\infty}$ norm 
appears to reach an asymptotic value corresponding to a slightly 
faster soliton. This is due to the fact that we considered a 
perturbation of almost 10\%. The final state is thus a solitary wave 
with larger mass, the solitary wave appears to be stable. 

\begin{figure}[htb!]
 \includegraphics[width=0.49\textwidth]{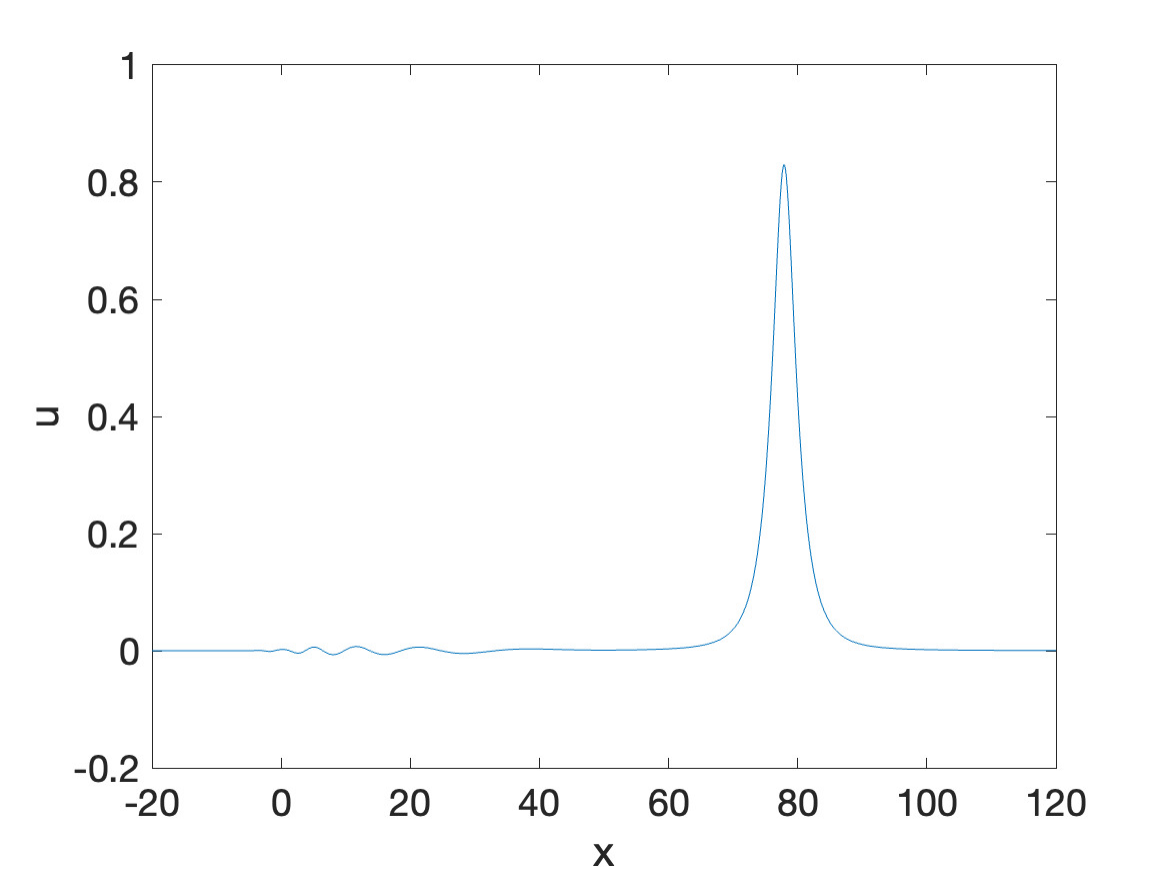}
 \includegraphics[width=0.49\textwidth]{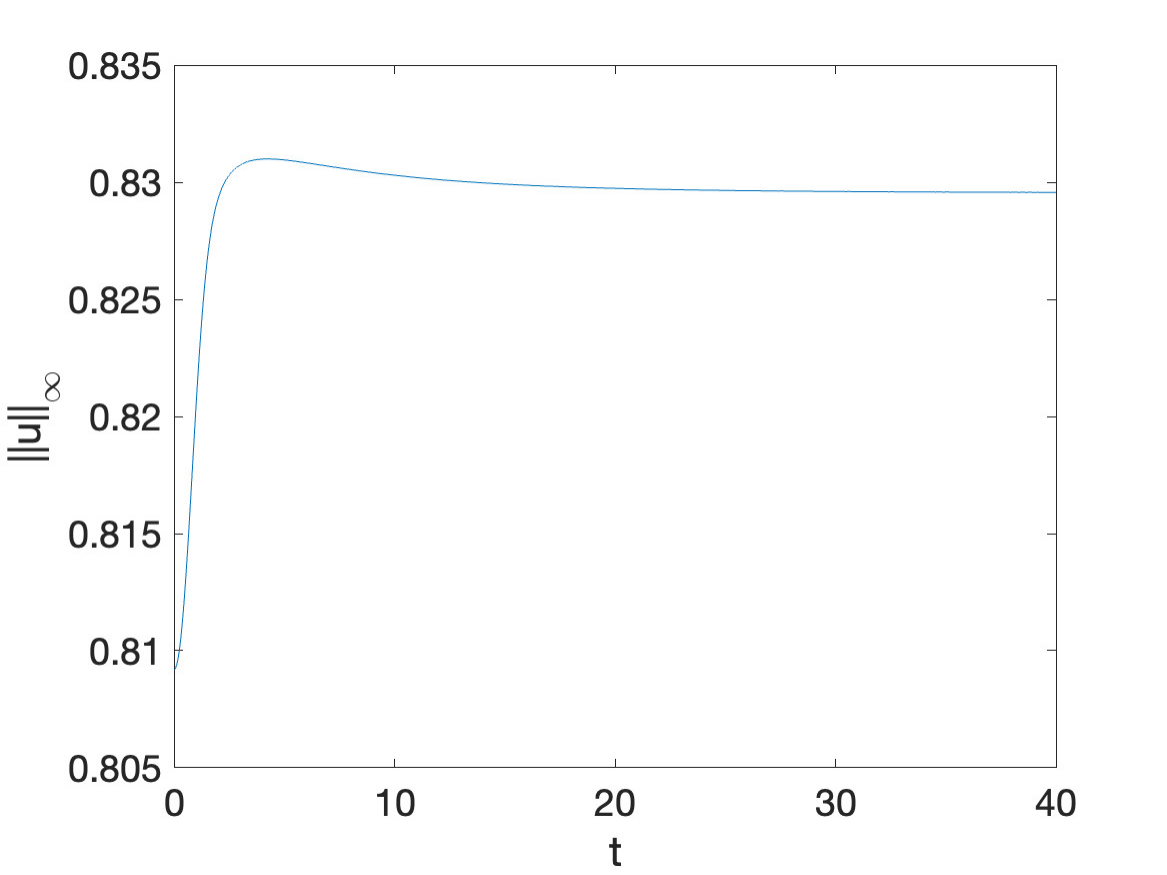}
 \caption{Solution to the fractional CH equation for initial data of 
 the form (\ref{gausspert}) with $A=-0.08$ for $\alpha=1.5$ and 
 $c=2$: on the left the solution for $t=40$, on the right the 
 $L^{\infty}$ norm of the solution in dependence of time. }
 \label{Qa15mgausspert}
\end{figure}

We obtain a similar result if we consider a perturbation of the 
solitary wave with slightly smaller mass than the unperturbed 
solitary wave. In Fig.~\ref{Qa15mgausspert}, we choose initial data 
of the form (\ref{gausspert}) with $A=-0.08$ with the same parameters 
as in Fig.~\ref{Qa15gausspert}. The solution at the final time $t=40$ 
appears to be again a solitary wave plus radiation. This is confirmed 
by the $L^{\infty}$ norm on the right of the same figure. The final 
state of the solution appears to be a solitary wave with a slightly 
smaller mass and velocity than $Q_{2}$. The solitary wave seems to 
be again stable even against comparatively large perturbations.

\begin{figure}[htb!]
 \includegraphics[width=0.49\textwidth]{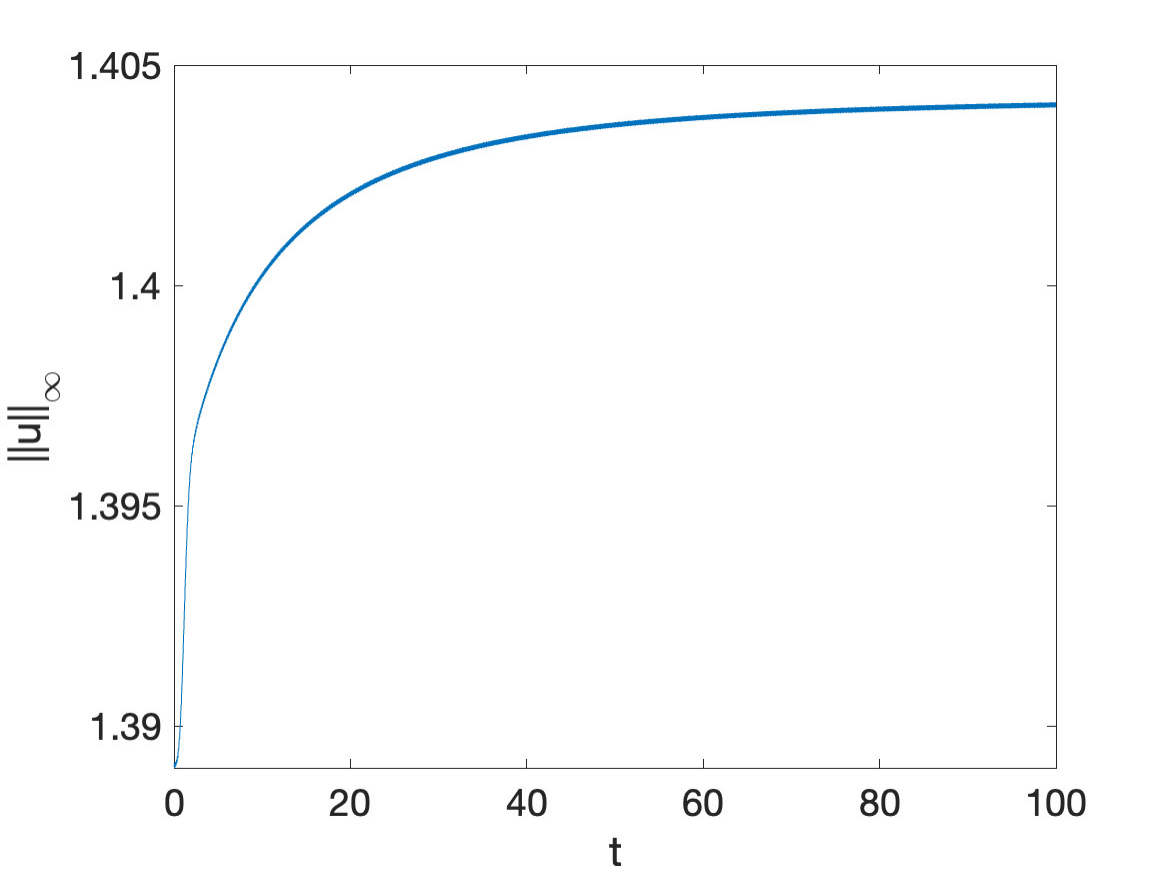}
 \includegraphics[width=0.49\textwidth]{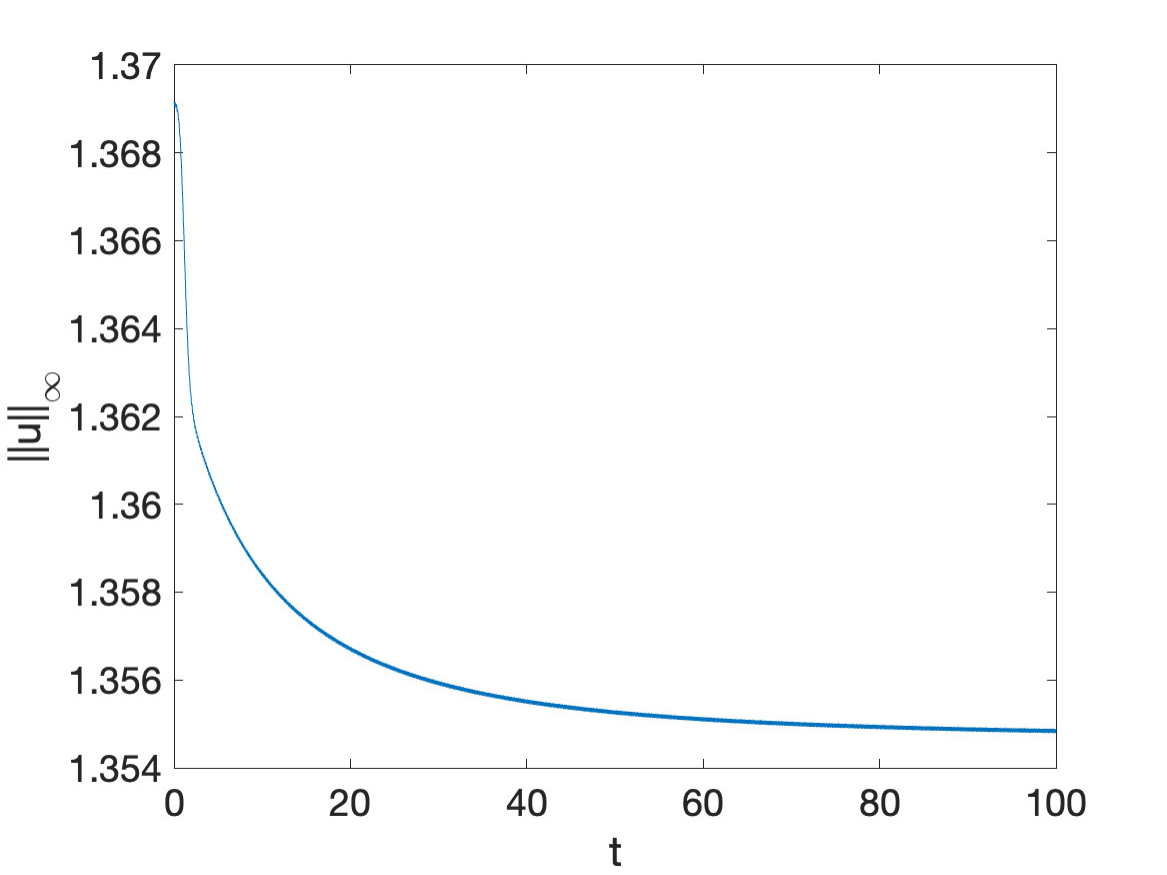}
 \caption{$L^{\infty}$ norm of the solution to the fractional 
 CH equation for initial data of 
 the form (\ref{gausspert}) for $\alpha=0.9$ and 
 $c=2$: on the left for $A=0.01$, on the right for $A=-0.01$. }
 \label{Qa09gausspert}
\end{figure}

It is an interesting question whether a similar behavior can be 
observed for smaller values of $\alpha$, i.e., for a fractional CH 
equation with less dispersion. We consider the case $\alpha=0.9$ for 
which we could construct solitary waves with $c=2$. Here we consider 
smaller perturbations of the order of 1\% for larger times than 
before. We apply $N_{t}=2*10^{4}$ time steps for $t\leq100$ to 
initial data of the form (\ref{gausspert}). In 
Fig.~\ref{Qa09gausspert} we show on the left the $L^{\infty}$ norm of 
the solution for $A=0.01$ and on the right for $A=-0.01$. In both cases 
the final state of the solution appears to be a solitary wave plus 
radiation. The small oscillations in the $L^{\infty}$ norm are due to 
the fact that we are working on a torus where the radiation can 
reenter the computational domain on the other side and that we 
determine the maximum of the solution on grid points. 

\begin{figure}[htb!]
 \includegraphics[width=0.49\textwidth]{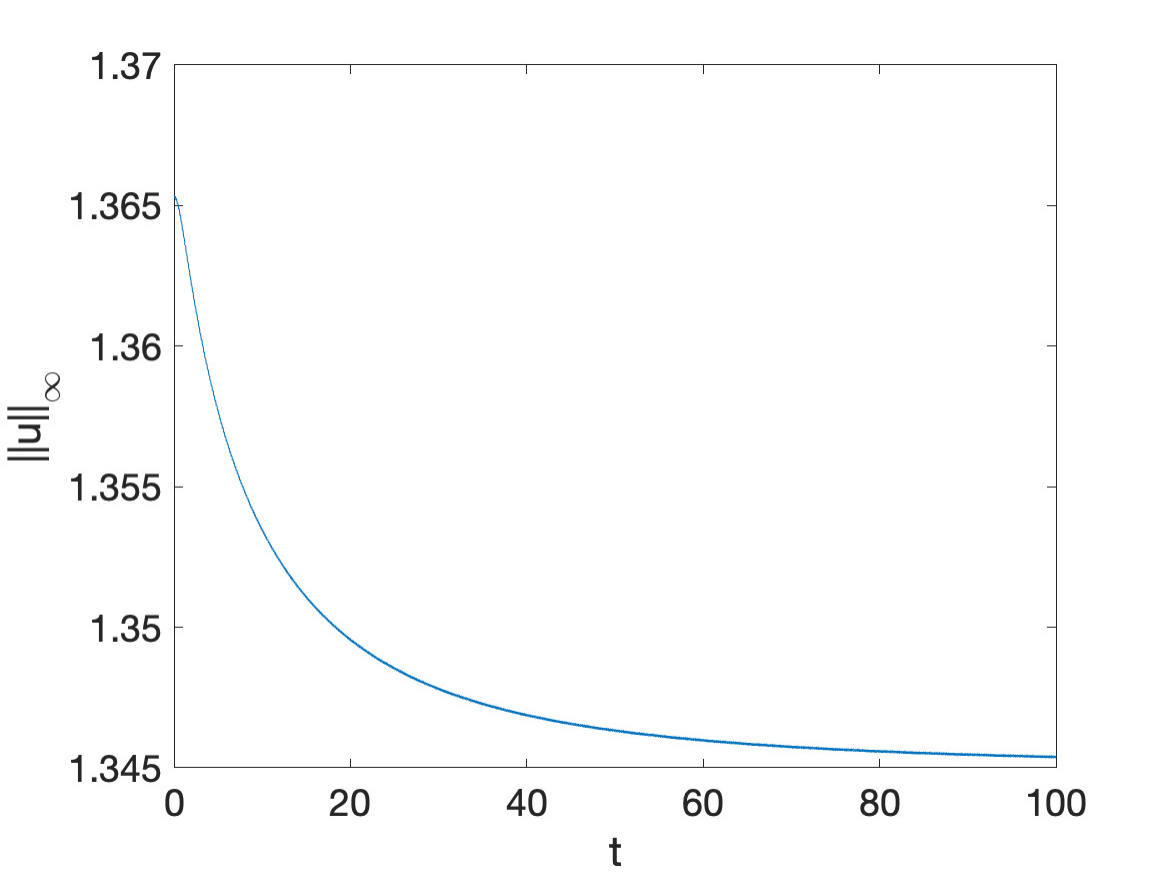}
 \includegraphics[width=0.49\textwidth]{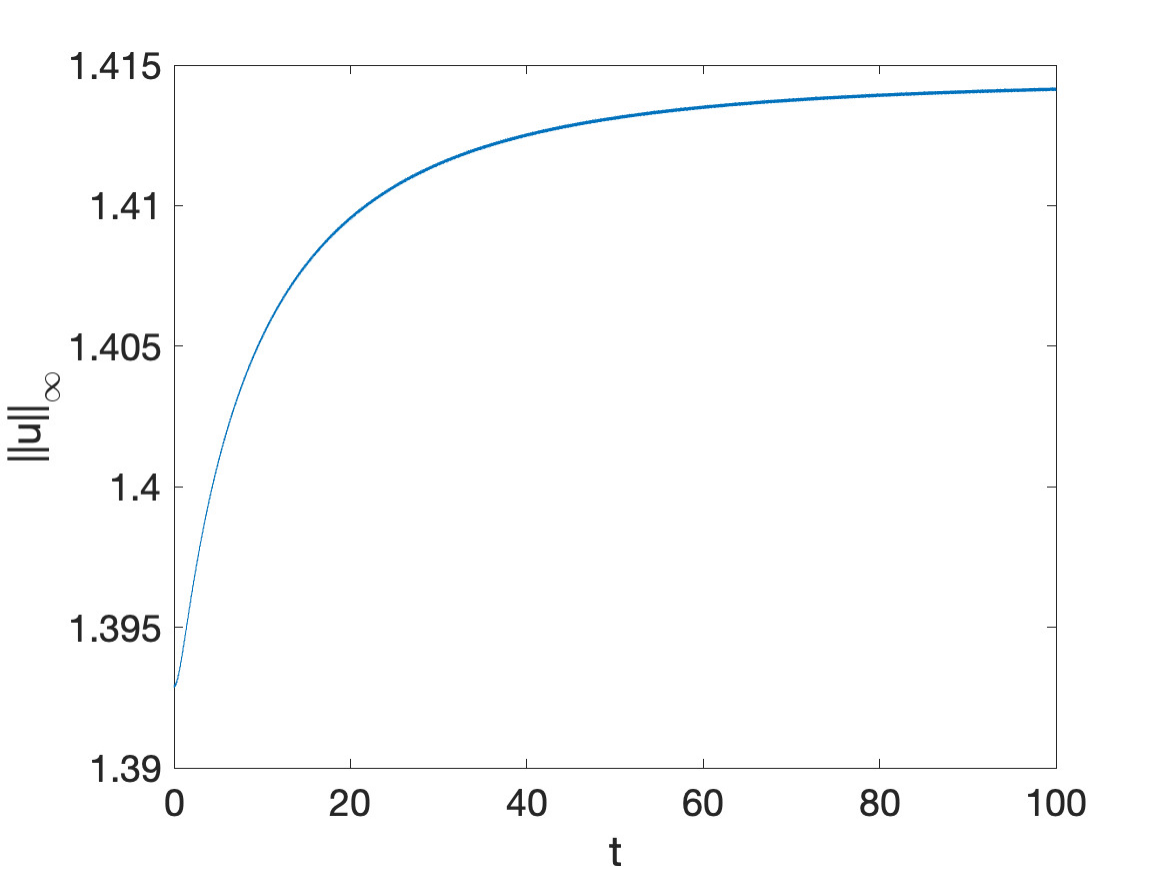}
 \caption{$L^{\infty}$ norm of the solution to the fractional 
 CH equation for initial data of 
 the form $u(x,0)=\lambda Q_{2}(x)$ for $\alpha=0.9$ and 
 $c=2$: on the left for $\lambda=0.99$, on the right for $\lambda=1.01$. }
 \label{Qa09lapert}
\end{figure}

The picture is very similar for different types of perturbations.  We 
consider the solution for intial data of the 
form $u(x,0) = \lambda Q_{2}(x)$ for $\alpha=0.9$ and real $\lambda$ 
close to 1. In 
Fig.~\ref{Qa09lapert}, the $L^{\infty}$ norms of the solution in 
dependence of time is shown, on the left for $\lambda=0.99$, on the 
right for $\lambda=1.01$. In both cases the final state appears to be 
a solitary wave of slightly different mass. 

\section{Localized initial data}
In this section we study the long time behavior of initial data from 
the Schwartz class of smooth rapidly decreasing functions. We are interested whether solitary waves 
appear in the solution asymptotically as expected from the \emph{soliton 
resolution conjecture}, or whether there is a \emph{blow-up}, a loss 
of regularity of the solution in finite time. 

\begin{figure}[htb!]
 \includegraphics[width=\textwidth]{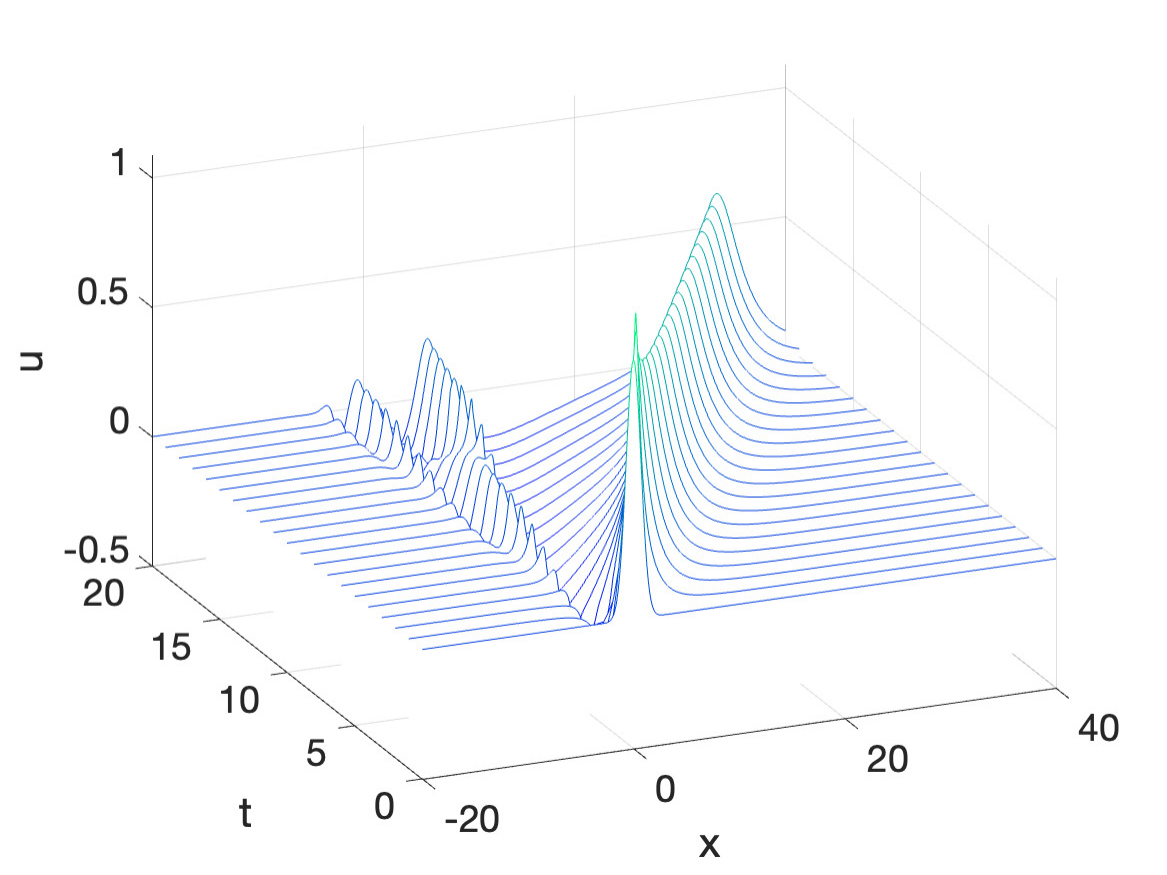}
 \caption{Solution to the fractional CH equation with $\alpha=1.5$ for initial data 
 $u(x,0)=\exp(-x^{2})$. }
 \label{fCHa15gausswater}
\end{figure}

Concretely we will study
Gaussian initial data, 
\begin{equation}
	u(x,0) = A\exp(-x^{2})
	\label{gauss},
\end{equation}
where $A>0$ is constant. We will use $N_{t}=10^{4}$ time steps and 
$N=2^{14}$ Fourier modes for  $x\in 3[-\pi,\pi]$. 
We first study the 
case $\alpha=1.5$. Small initial data will be simply radiated towards 
infinity. But if we take initial data of the form (\ref{gauss}) with 
$A=1$, we get the solution shown in Fig.~\ref{fCHa15gausswater}. The 
initial hump breaks up into several humps. One of them, possibly 
three, appear to be solitary waves.  It seems that the solution indeed 
decomposes into solitary waves and radiation. The precise number of  
solitary waves appearing for large times is unknown. 

\begin{figure}[htb!]
 \includegraphics[width=0.49\textwidth]{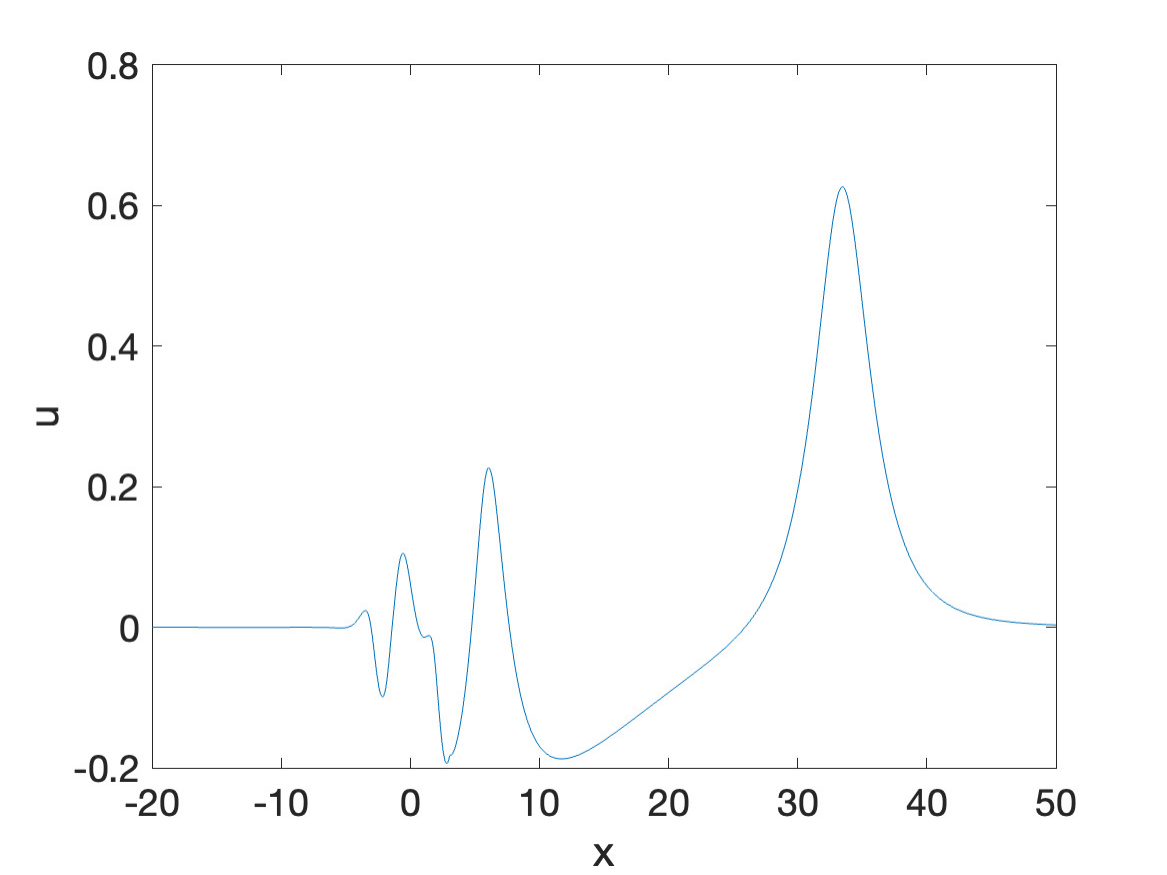}
 \includegraphics[width=0.49\textwidth]{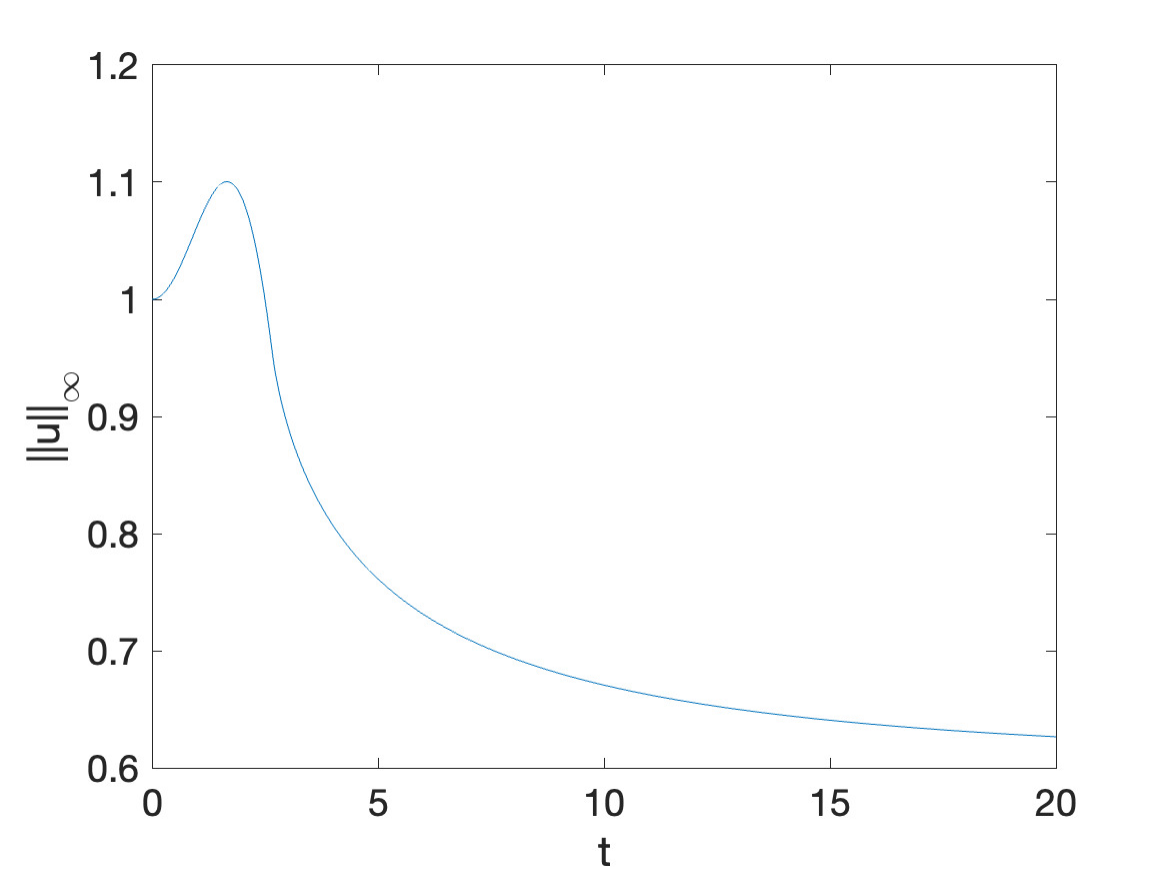}
\caption{Solution to the fractional CH equation with $\alpha=1.5$ for initial data 
 $u(x,0)=\exp(-x^{2})$, on the right the $L^{\infty}$ norm, on the 
 left the solution for $t=20$. }
 \label{fCHa15gauss}
\end{figure}

The interpretation of the largest hump as a solitary wave is 
confirmed by the $L^{\infty}$ norm of the solution shown in 
Fig.~\ref{fCHa15gauss} on the right. It seems to reach a constant 
asymptotic value as 
expected for a solitary wave. Since the solitary waves for the 
fractional CH equation do not have a simple scaling with the velocity 
$c$ as the fKdV solitary wave, a fit of the hump to the solitary waves is not 
obvious. The solution for $t=20$ is shown on the left of 
Fig.~\ref{fCHa15gauss}. 

\begin{figure}[htb!]
 \includegraphics[width=0.49\textwidth]{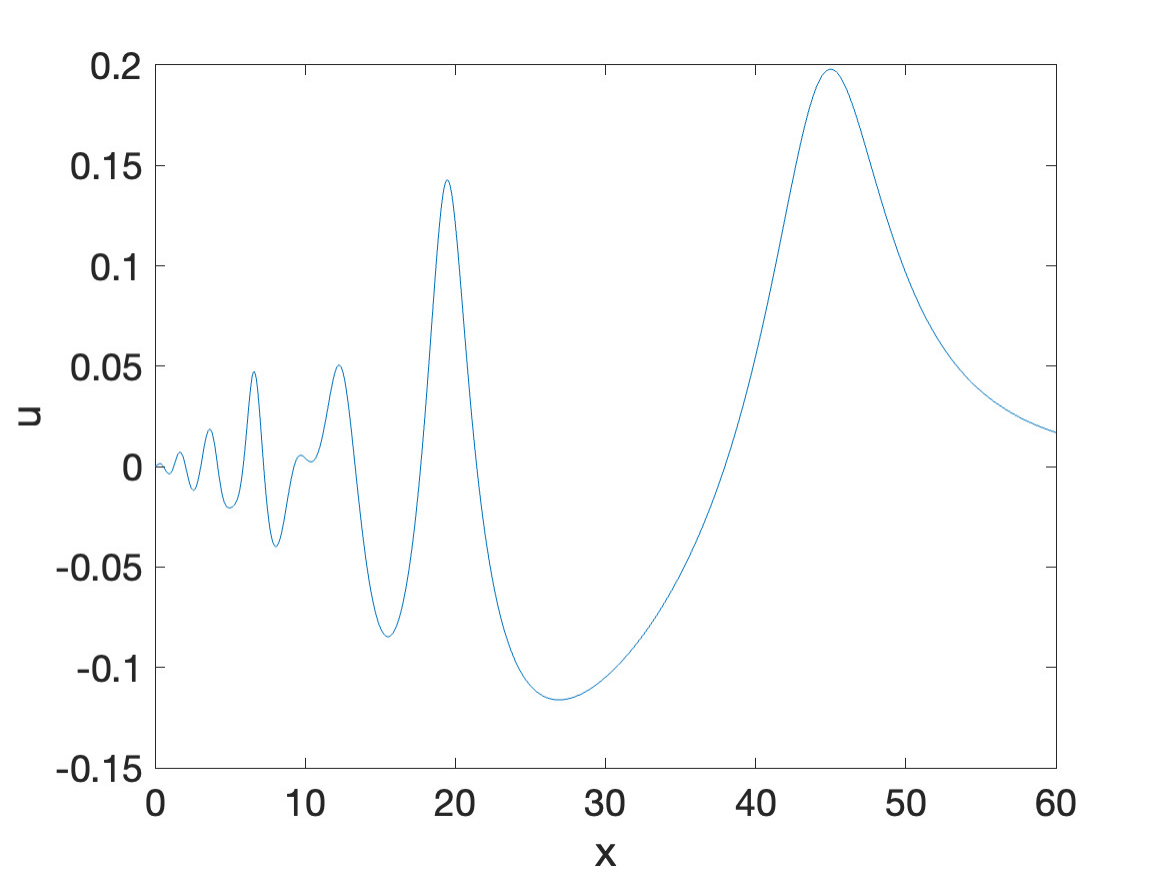}
 \includegraphics[width=0.49\textwidth]{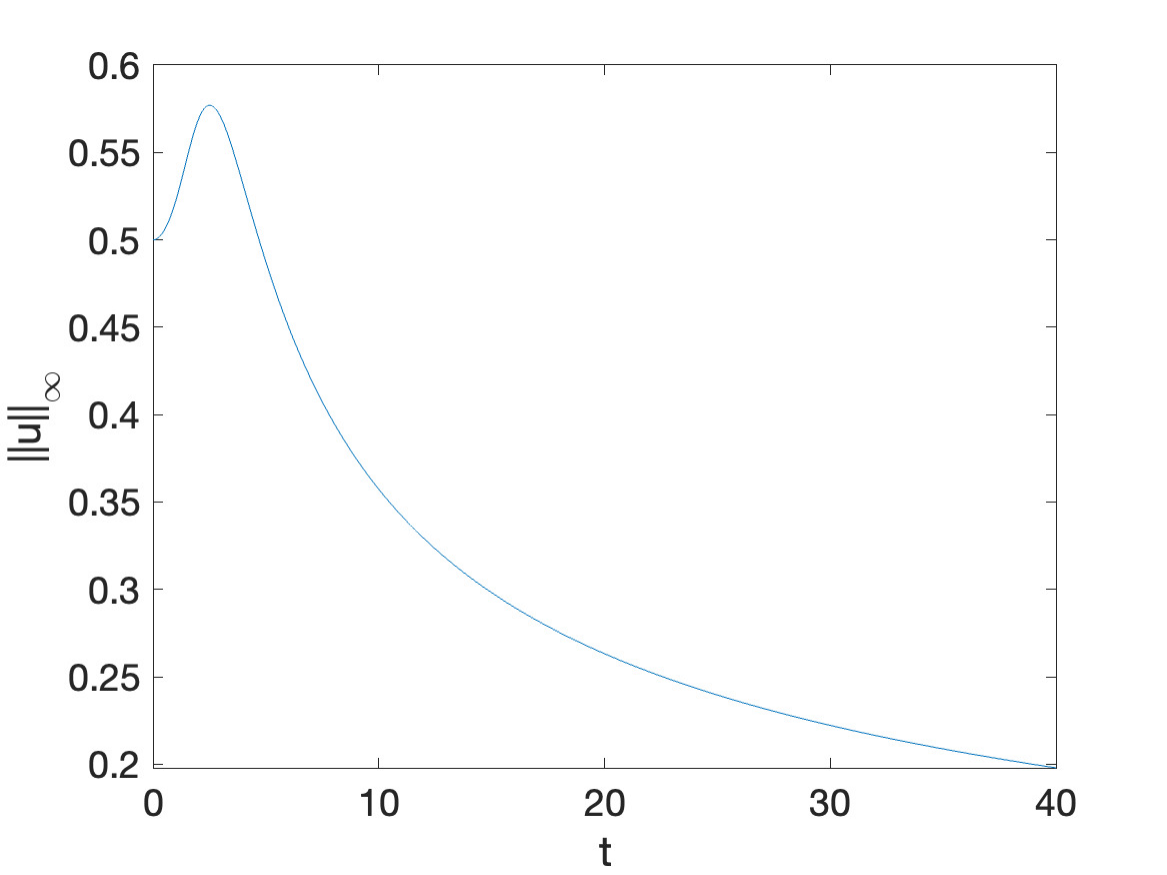}
 \caption{Solution to the fractional CH equation with $\alpha=0.9$ for initial data 
 $u(x,0)=0.5\exp(-x^{2})$, on the left the solution for $t=40$, on the 
 right the $L^{\infty}$ norm. }
 \label{fCHa0905gauss}
\end{figure}

For smaller values of $\alpha$, the picture can change. If we take 
$\alpha=0.9$ and initial data of the form (\ref{gauss}) with $A=0.5$, 
the initial hump will be radiated away. The solution for $t=40$ is 
shown on the left of Fig.~\ref{fCHa0905gauss}. The $L^{\infty}$ norm of 
the solution on the right of the same figure appears to decrease 
monotonically with time. 

\begin{figure}[htb!]
 \includegraphics[width=0.49\textwidth]{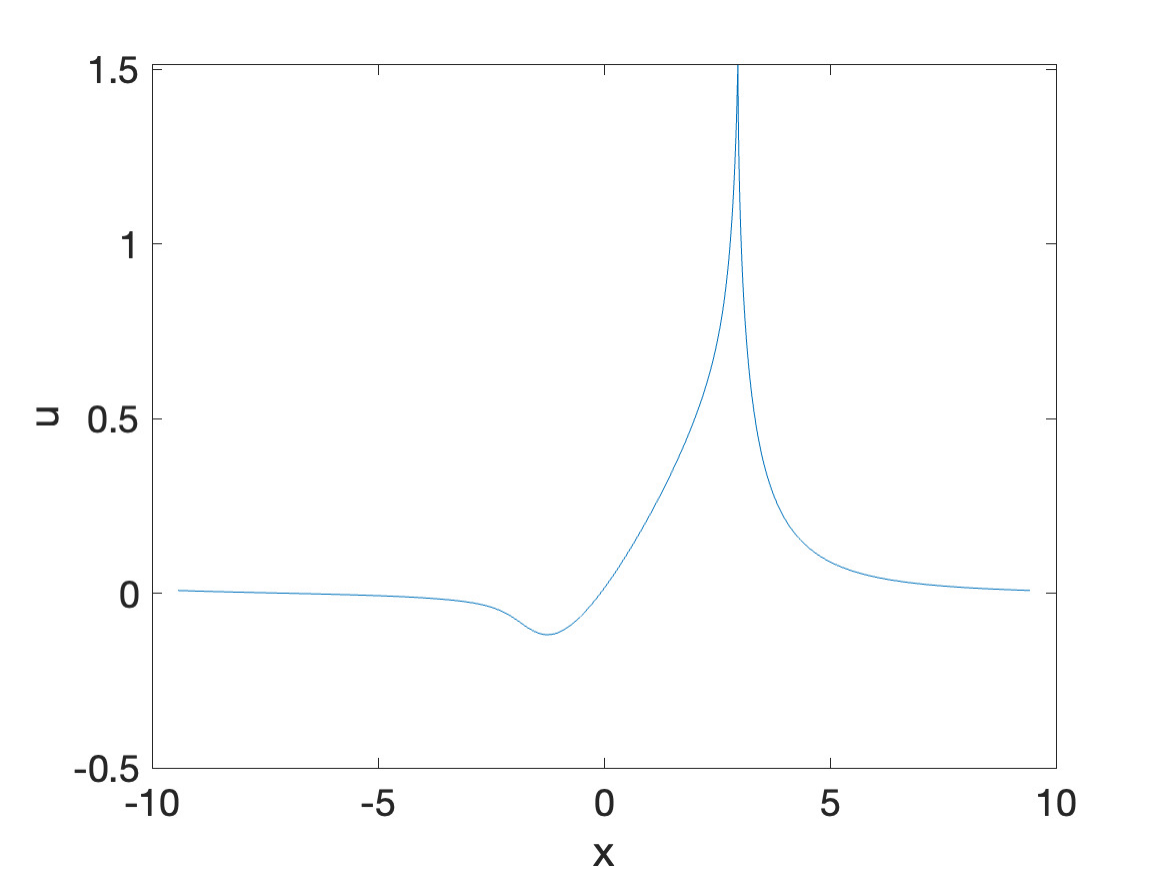}
 \includegraphics[width=0.49\textwidth]{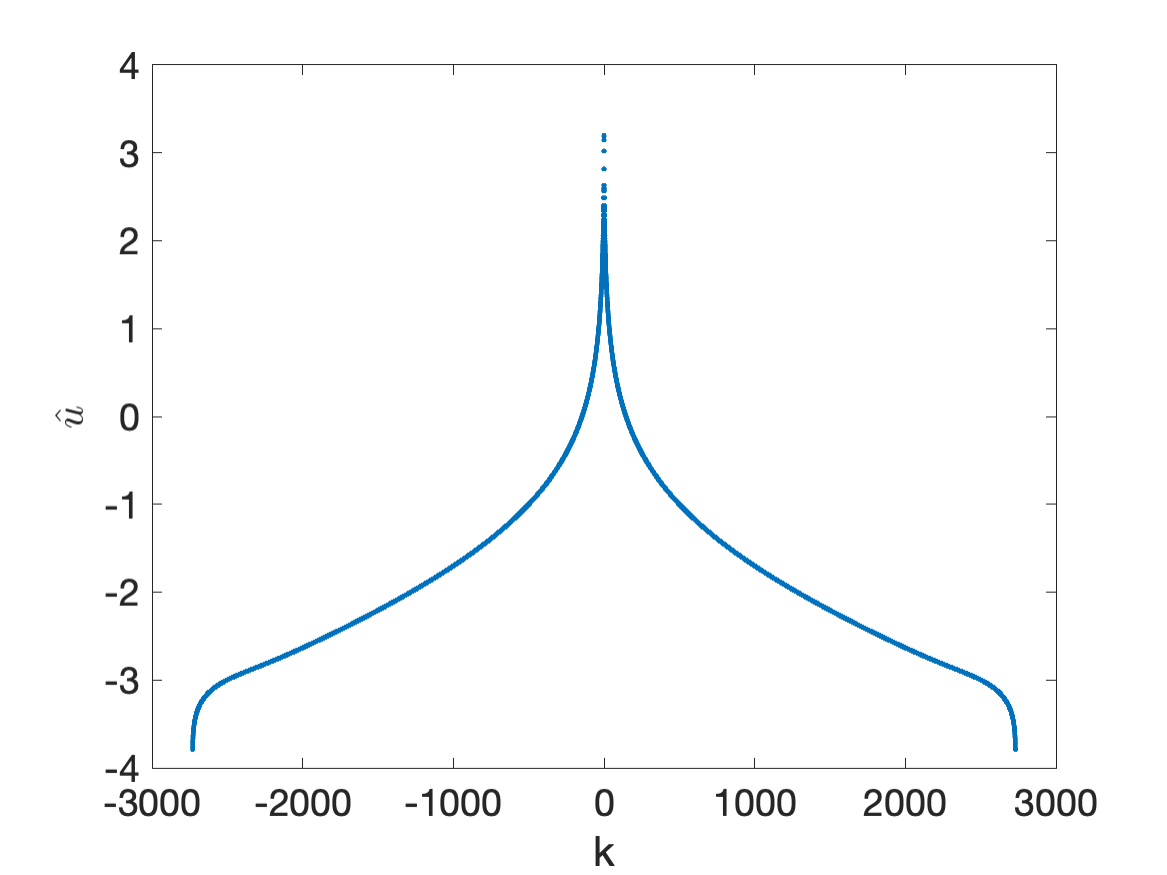}
 \caption{Solution to the fractional CH equation with $\alpha=0.9$ for initial data 
 $u(x,0)=\exp(-x^{2})$, on the left the solution for $t=1.7667$, on the 
 right the DFT coefficients for this solution. }
 \label{fCHa09gauss}
\end{figure}

However for initial data of the form (\ref{gauss}) of larger mass, i.e., 
larger values of $A$, we do not find solitary waves in the long time 
behavior of the solution. Instead for $A=1$, a cusp appears to form 
in finite time, for $t\sim 1.7667$ in this case. The solution at this 
time can be seen in Fig.~\ref{fCHa09gauss} on the left. 

It is not only the figure that indicates a cusp formation, this is also 
confirmed by the Fourier coefficients (more precisely the DFT of $u$) 
on the right of the same figure. The algebraic decay of these 
coefficients with the index $|k|$ indicates that a singularity of the 
analytic continuation of the function $u$ to the complex plane will 
hit at the critical time the real axis. It is well 
known that an essential singularity in the complex plane of the form $u\sim 
(z-z_{j})^{\mu_{j}}$, $\mu_{j}\notin \mathbb{Z}$, 
with $z_{j}=\alpha_{j}-i\delta_{j}$ in the lower half 
plane ($\delta_{j}\geq 0$) implies    for $k\to\infty$  
the following asymptotic behavior of the Fourier 
coefficients (see e.g.~\cite{asymbook}),
\begin{equation}
    \hat{u}\sim 
    \sqrt{2\pi}\mu_{j}^{\mu_{j}+\frac{1}{2}}e^{-\mu_{j}}\frac{(-i)^{\mu_{j}+1}}{k^{\mu_{j}+1}} e^{-ik\alpha_{j}-k\delta_{j}}
    \label{fourierasym}.
\end{equation}
For a single such singularity with positive $\delta_{j}$, the modulus of the Fourier 
coefficients decreases exponentially for large $k$ until
$\delta_{j}=0$, when 
the modulus of the Fourier coefficients has an algebraic dependence 
on $k$. As first shown in \cite{SSF}, one can use (\ref{fourierasym}) 
to numerically characterize the singularity as discussed in detail in 
\cite{KR} (the DFT coefficients have a similar behavior as the 
Fourier transform). Essentially a least square method is applied to 
$\ln|\hat{u}|$ to fit the parameters for $k>100$. The reader is 
referred to \cite{KR} for details. For the case shown in 
Fig.~\ref{fCHa09gauss}, we find $\mu=0.5007$. This means that a 
square root type singularity is observed in this case which provides 
numerical evidence of the second part of the main conjecture in the 
introduction. Since  global existence in time does not appear to hold 
for solutions to initial 
data of sufficiently lage mass, no solitary waves are observed in 
this case. 

Note that this is in apparent contradiction to the stability of the 
solitary waves observed numerically in the previous section. However, 
it has to be remembered that the solitary waves have a slow 
algebraic fall-off towards infinity whereas we consider 
exponentially decaying data in this section. For $L^{2}$-critical 
generalized 
Korteweg-de Vries equations, a blow-up in finite time is only 
observed if the intitial data are sufficiently rapidly decreasing, 
see \cite{MMR} and references therein. The situation appears to be 
similar for fCH solutions for sufficiently small $\alpha$. 

\section{Dispersive shock waves}
In this section we study the appearence of dispersive shock waves 
(DSWs) in 
fCH solutions. A convenient way to study the formation of zones of 
rapid oscillations in the solutions of dispersive PDEs is to consider 
the solution for large times on large scales. This can be done by 
introducing a small parameter $\epsilon\ll1$ and rescale $t$, $x$
according to $t\mapsto t/\epsilon$, $x\mapsto x/\epsilon$. This leads 
for equation (\ref{fCH}) to
\begin{equation}\label{fCHe}
u_t+ \kappa_1 u_x + 3 uu_x + \epsilon^{\alpha}D^{\alpha} u_t = 
-\kappa_2 \epsilon^{\alpha}[2D^{\alpha} 
(u u_x)+ uD^{\alpha}u_x ],
\end{equation}
where we have kept the same notation as for the case $\epsilon=1$. 
Thus equation (\ref{fCHe}) is simply equation (\ref{fCH}) with 
$D^{\alpha}$ replaced by $\epsilon^{\alpha} D^{\alpha}$. In the 
formal limit $\epsilon\to0$, equation (\ref{fCHe}) reduces to the 
Hopf equation 
\begin{equation}
	u_t+ \kappa_1 u_x + 3 uu_x  = 0,
	\label{hopf}
\end{equation}
where the term linear in $u_{x}$ can be absorbed by a Galilei 
transformation. The Hopf equation is known to develop a gradient 
catastrophe for hump-like initial data in finite time. Dispersive 
regularisations of this equation as (\ref{fCHe}) are expected to lead 
to solutions with zones of rapid oscillations in the vicinity of the 
shocks of the Hopf solution for the same initial data. 

In \cite{GK}, we have studied numerically the onset of oscillations 
in solutions of the CH equation. It was conjectured that a special 
solution to the second equation in the Painlev\'e I hierarchy, see 
for instance \cite{GKK}, 
provides an asymptotic description of the break-up of CH solutions in 
this case. For larger times the oscillatory zone was studied in 
\cite{AGK}. A conjecture to describe the leading edge of the 
oscillatory zone in terms of a Painlev\'e transcendent was given. 

\begin{figure}[htb!]
 \includegraphics[width=0.49\textwidth]{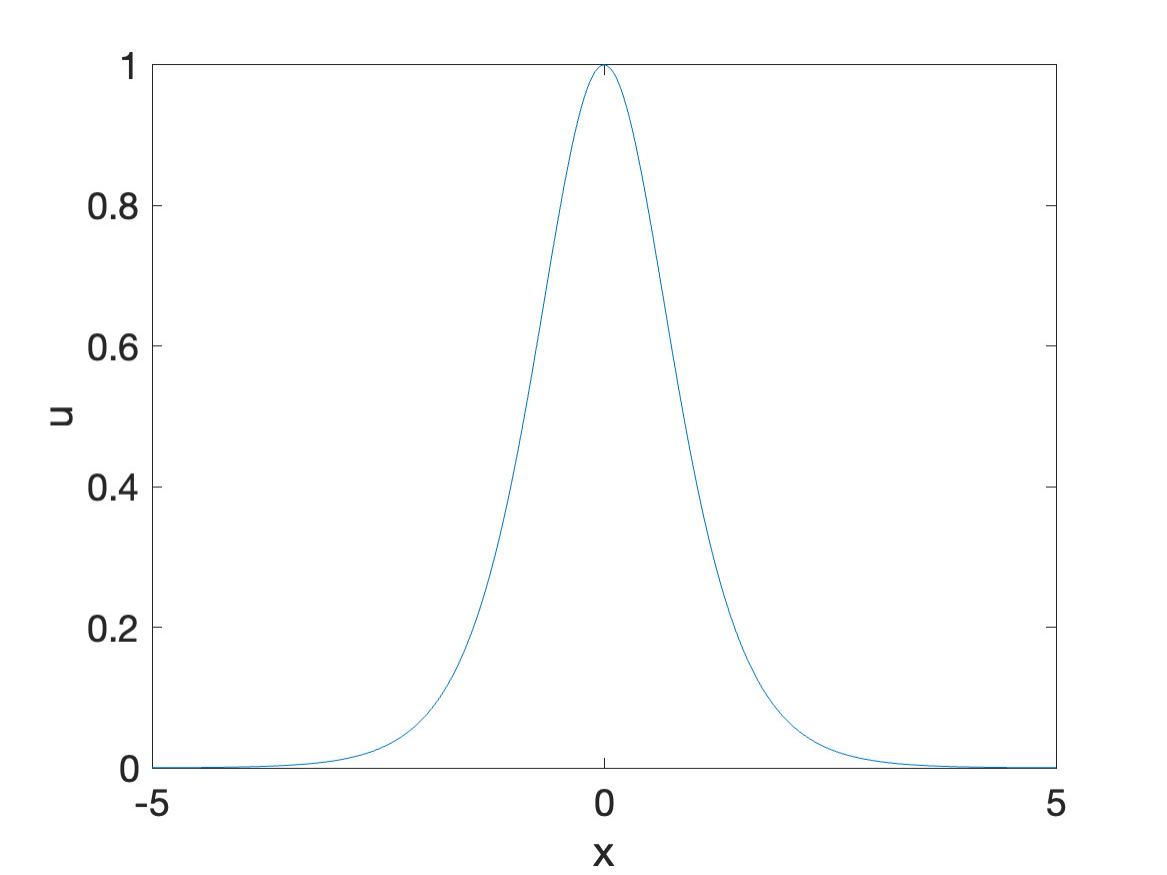}
 \includegraphics[width=0.49\textwidth]{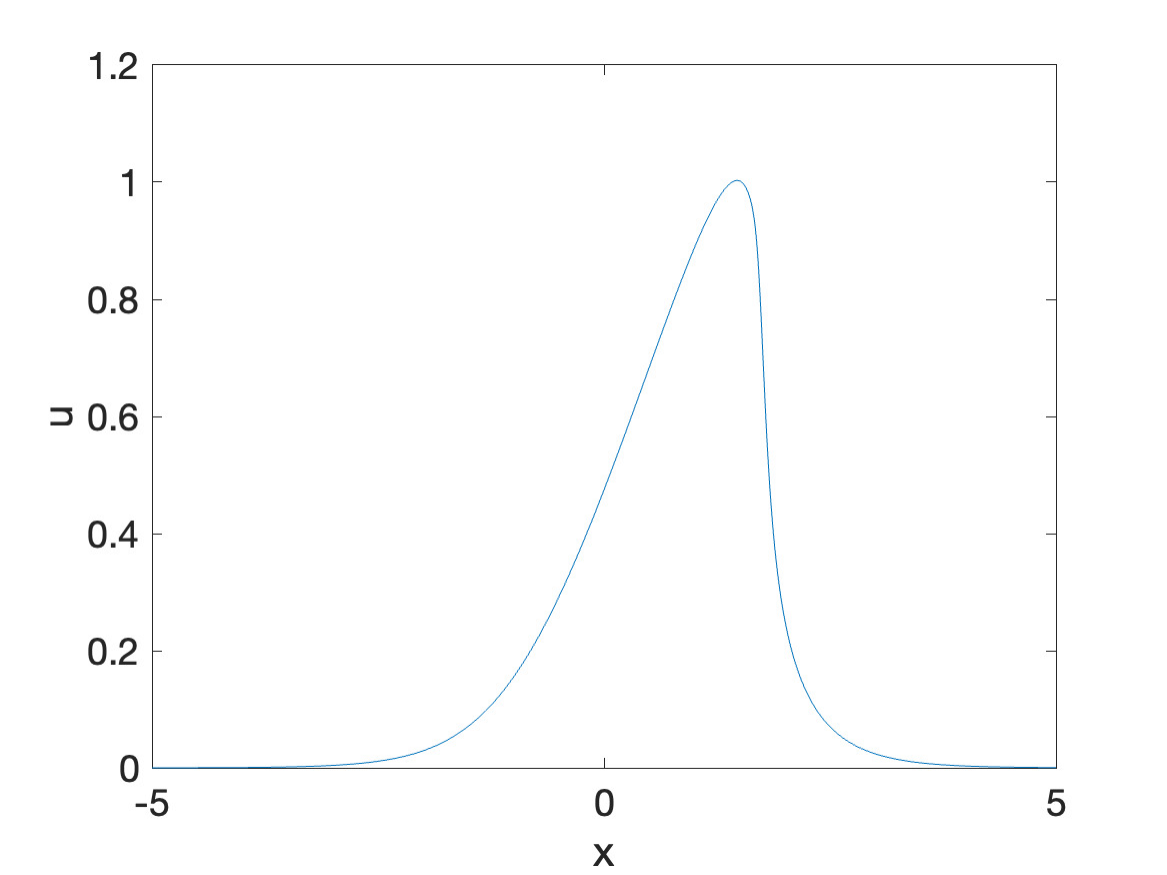}\\
  \includegraphics[width=0.49\textwidth]{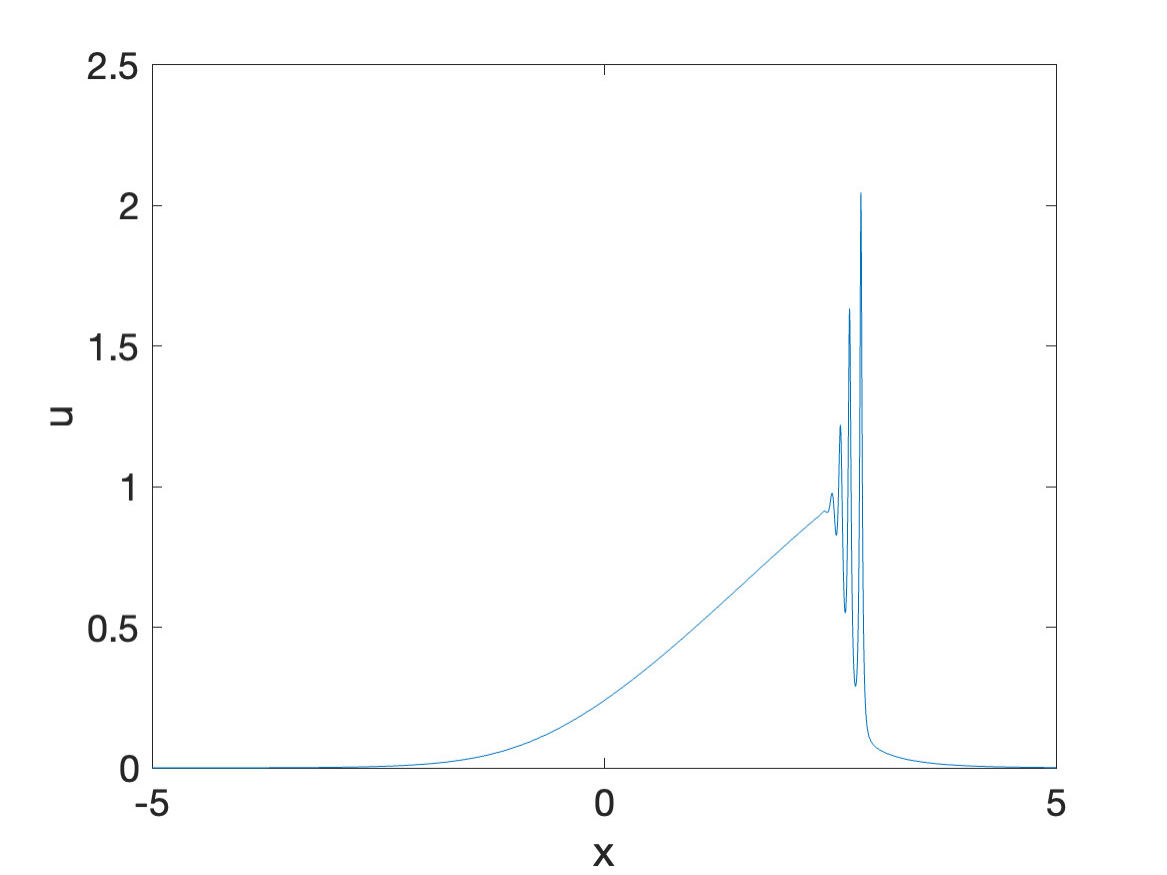}
 \includegraphics[width=0.49\textwidth]{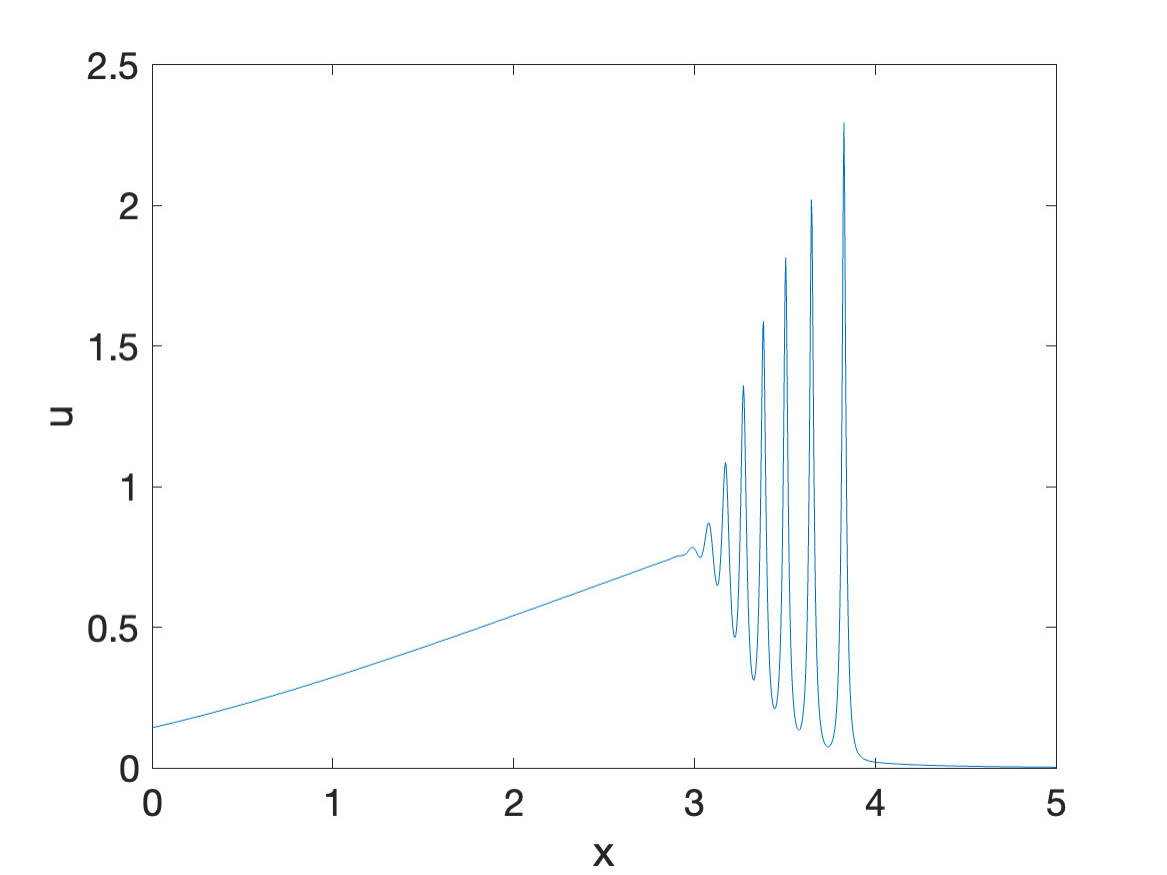}
 \caption{Solution to the fractional CH equation (\ref{fCHe}) with 
 $\alpha=1.5$, $\omega=0.6$ and $\epsilon=10^{-2}$ for initial data 
 $u(x,0)=\mbox{sech}^{2} x$ for several values of time (on the top 
 for $t=0$, $t=0.35$, in the bottom for $t=0.7$, $t=1$). }
 \label{fChsech1e2}
\end{figure}

We will study below similar examples as in \cite{GK,AGK} for fCH, initial 
data of the form $u(x,0)=\mbox{sech}^{2}x$ for several values of 
$\epsilon$. In Fig.~\ref{fChsech1e2}, we show the fCH solution for 
$\alpha=1.5$ and $\epsilon=10^{-2}$ for several values of $t$. We use 
$N=2^{14}$ Fourier modes and $N_{t}=10$ time steps for $t\leq1$. 
A first oscillation forms for $t\sim 0.4$ (the critical time for the 
Hopf solution is $t_{c}\sim 0.433$), then a well defined zone 
of oscillations appears also called the Whitham zone. 

\begin{figure}[htb!]
 \includegraphics[width=0.32\textwidth]{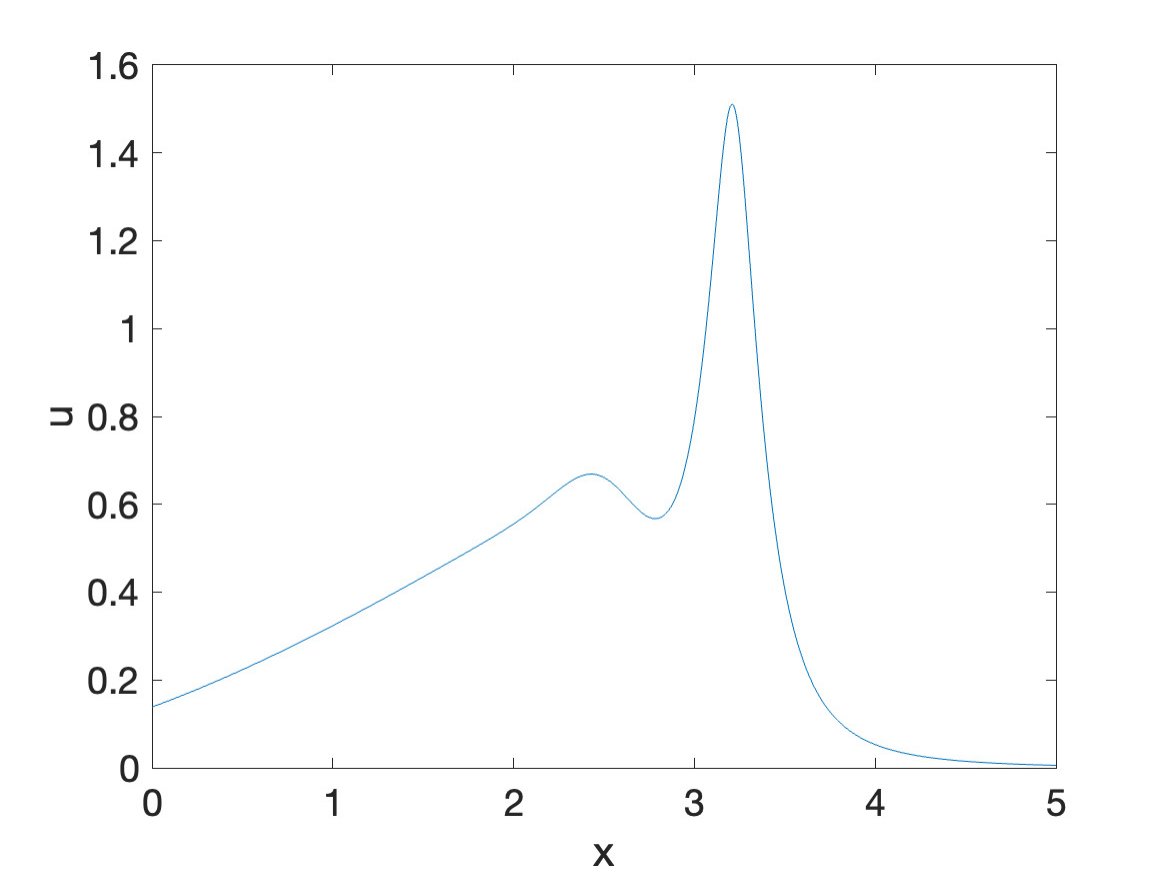}
  \includegraphics[width=0.32\textwidth]{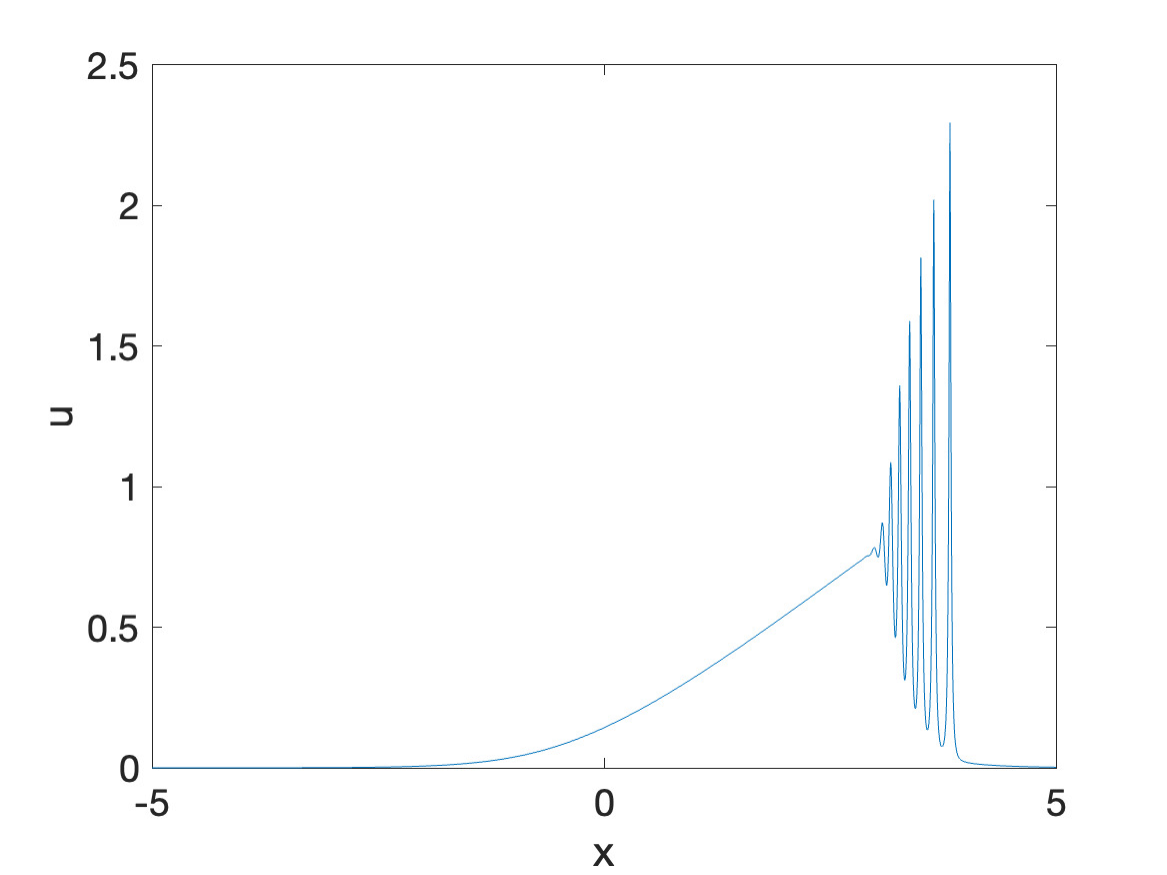}
 \includegraphics[width=0.32\textwidth]{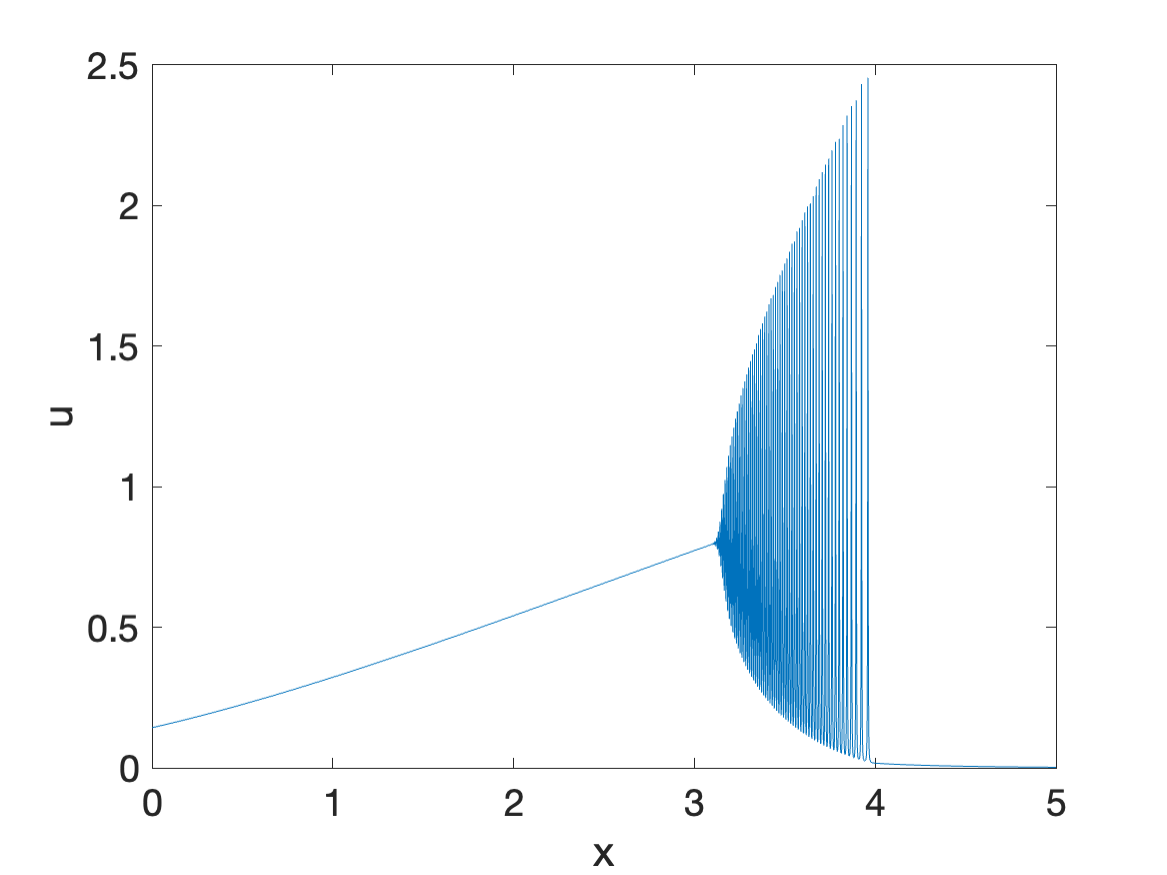}
 \caption{Solution to the fractional CH equation (\ref{fCHe}) with 
 $\alpha=1.5$, $\omega=0.6$  for initial data 
 $u(x,0)=\mbox{sech}^{2} x$ for  $t=1$ for $\epsilon=10^{-1}$, 
 $\epsilon=10^{-2}$, $\epsilon=10^{-3}$ from left to right. }
 \label{fChseche}
\end{figure}

The Whitham zone becomes more defined and more oscillatory the 
smaller $\epsilon$ is. Thus there is no strong limit $\epsilon\to0$ 
for DSWs. Note that DSWs in CH solutions are numerically more demanding than the 
corresponding KdV solutions for which special integrators exist, see 
for instance \cite{etna}. For CH time integration  is more problematic since the 
dispersive terms are nonlinear in contrast to KdV. The reduced 
dispersion in CH compared to KdV because of the nonlocality 
($(1-\partial_{xx})u_{t}$ leads to less oscillations in CH than in 
similar KdV situations, but to stronger gradients). This is amplified 
in fCH solutions since the dispersion is smaller than in CH. To treat 
the case $\epsilon=10^{-3}$ for the situation shown in 
Fig.~\ref{fChsech1e2}, we apply $N=2^{18}$ Fourier modes and 
$N_{t}=10^{5}$ time steps. We show the fCH solution for three values 
of $\epsilon$ in Fig.~\ref{fChseche} for the same initial data at the 
same time. 

\begin{figure}[htb!]
 \includegraphics[width=0.49\textwidth]{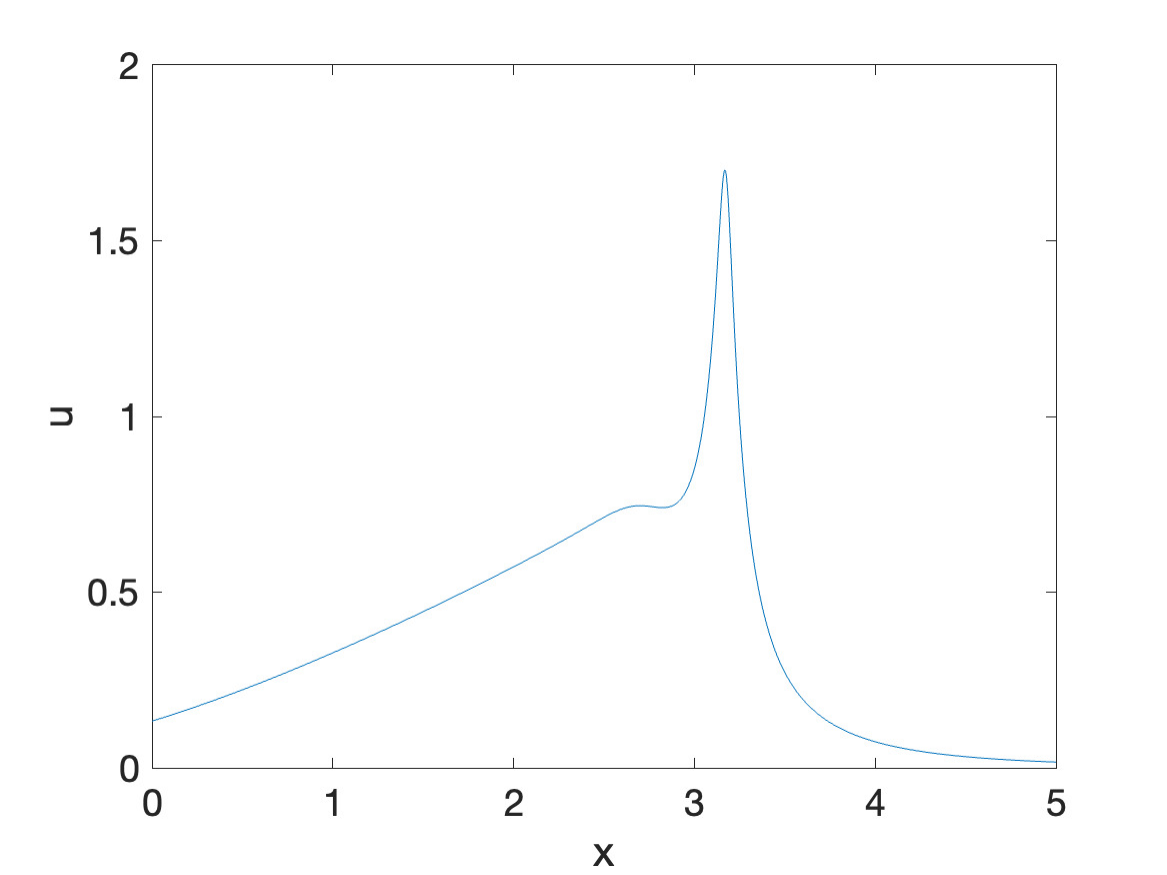}
 \includegraphics[width=0.49\textwidth]{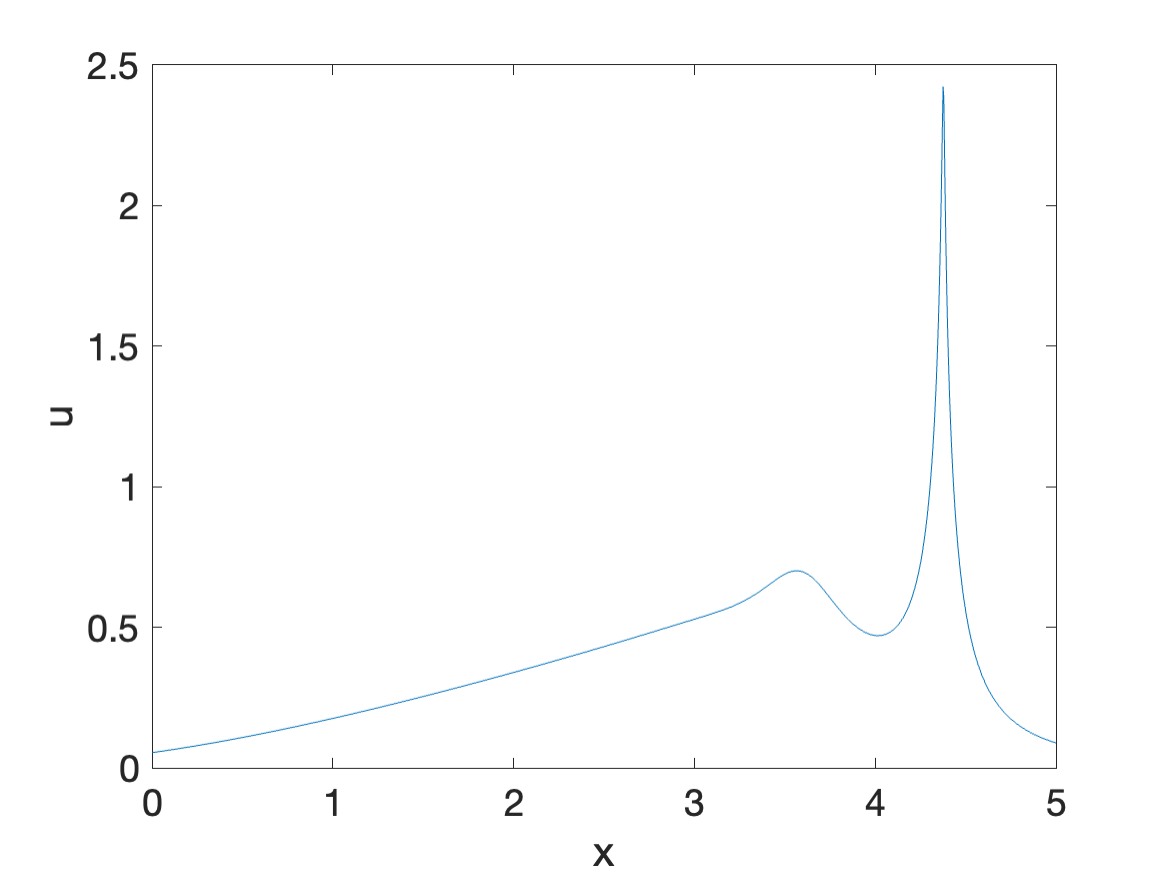}
 \caption{Solution to the fractional CH equation (\ref{fCHe}) with 
 $\alpha=0.9$, $\omega=0.6$ and $\epsilon=10^{-1}$ for initial data 
 $u(x,0)=\mbox{sech}^{2} x$ for  $t=1.005$ on the left and $t=1.5$ on 
 the right. }
 \label{fChsecha09e1}
\end{figure}

For smaller $\alpha$, the dispersion is even weaker. For the same 
initial data as in Fig.~\ref{fChseche}, we get for $\alpha=0.9$ 
and $\epsilon=10^{-1}$ again a DSW as can be seen in 
Fig.~\ref{fChsecha09e1}. As expected the smaller dispersion than in 
Fig.~\ref{fChseche} leads to less oscillatory behavior.

\begin{figure}[htb!]
 \includegraphics[width=0.49\textwidth]{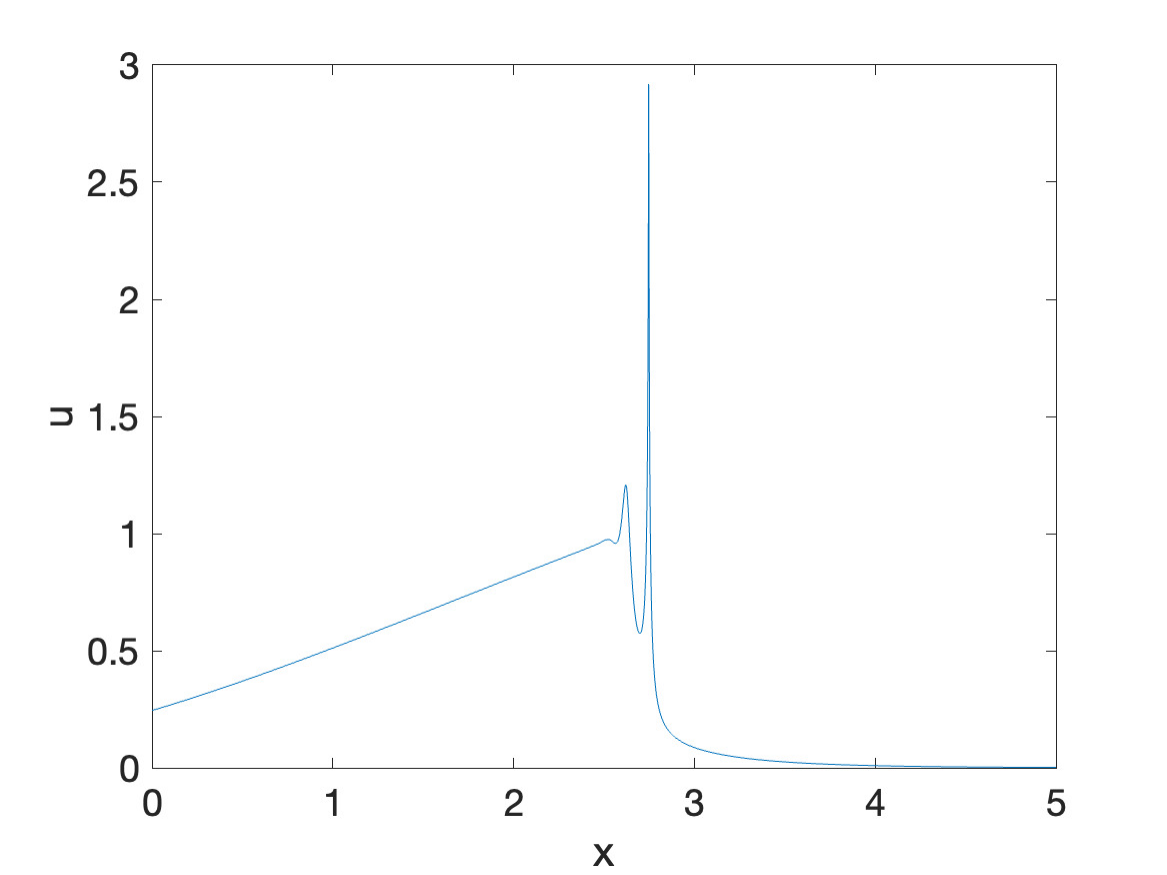}
 \includegraphics[width=0.49\textwidth]{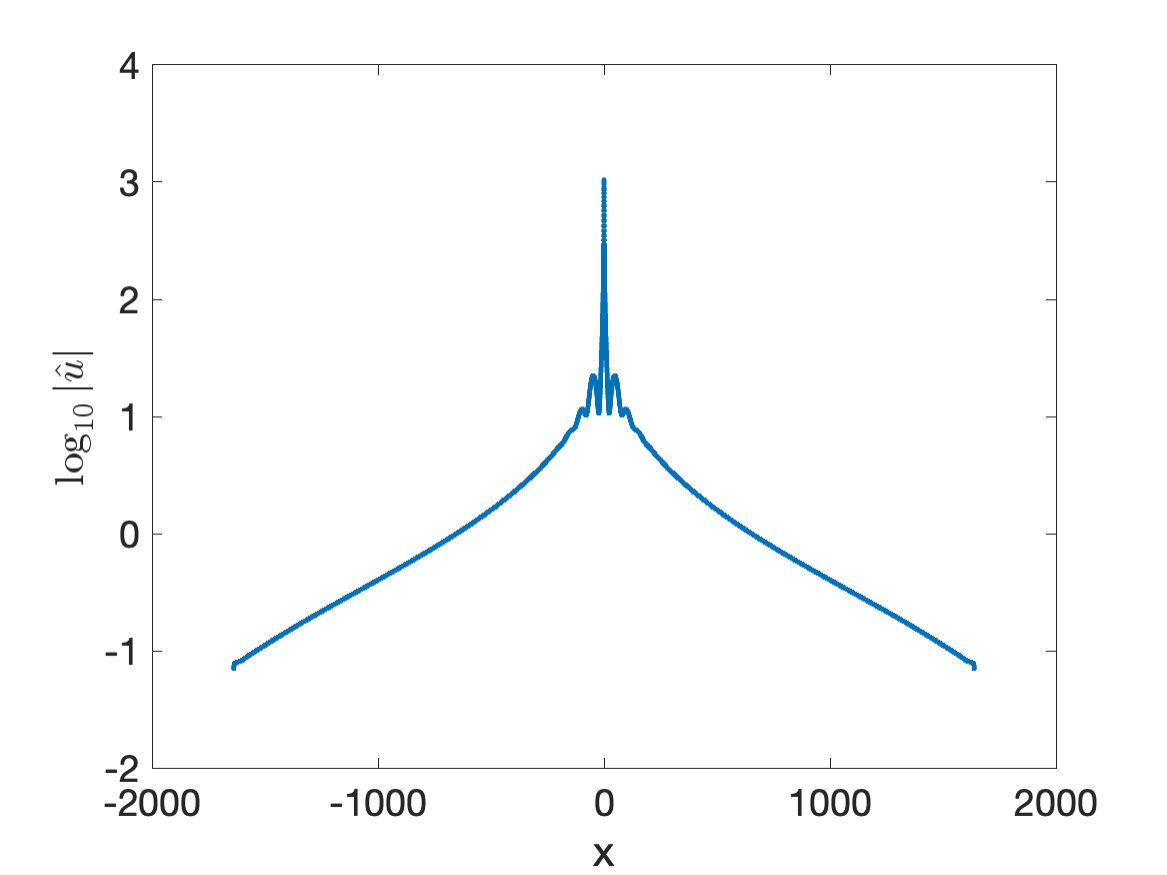}
 \caption{Solution to the fractional CH equation (\ref{fCHe}) with 
 $\alpha=0.9$, $\omega=0.6$ and $\epsilon=10^{-2}$ for initial data 
 $u(x,0)=\mbox{sech}^{2} x$ for  $t=0.6819$ on the left and the DFT 
 coefficients  on 
 the right. }
 \label{fChsecha09e2}
\end{figure}

However the situation is different for smaller $\epsilon$ than in 
Fig.~\ref{fChsecha09e1}, $\epsilon=10^{-2}$, as can be seen in 
Fig.~\ref{fChsecha09e2}. We use $N=2^{15}$ Fourier modes and 
$N_{t}=10^{4}$ time steps for $t\leq 0.7$. There is once more a DSW forming, but the 
first peak appears to develop in finite time into a cusp, see the 
left of Fig.~\ref{fChsecha09e2}. It is not surprising that a smaller 
$\epsilon$ leads to a cusp formation that was already observed in the 
previous section for initial data of sufficient mass in this case. 
Since the formal rescaling of $x$ and $t$ with $\epsilon$ leads also 
to a rescaling of the mass, the same initial data will have more mass 
in the original setting (\ref{fCH}) the smaller $\epsilon$. The cusp 
formation is confirmed by the DFT 
coefficients on the right of the figure which gives more numerical 
evidence to the second part of the Main Conjecture. A fitting of the 
coefficients to (\ref{fourierasym}) indicates an exponent $\mu\sim 
0.68$. This is slightly larger than the factor $1/2$ found in the 
previous section, but the accuracy in identifying the factor 
$\mu$ is always less than for the exponent $\delta$ in 
(\ref{fourierasym}), in particular here where the DSW already leads 
to a slower decay of the DFT coefficients with the index $k$. Thus 
there is strong evidence for cusp formation, but the exact character 
of the singularity needs to be justified analytically. 

\section{Outlook}
In this paper, we have  started a numerical study of the fractional 
CH equation. Solitary waves have been numerically constructed, and 
indications have been found that there could be a minimal value of 
$\alpha$ for given velocity $c$ and positive parameter $\omega$ below 
which there are no smooth solitary waves. It was shown that the 
numerically constructed smooth solitary waves are stable. 
A study of initial data from the Schwartz class of smooth rapidly 
decreasing functions led for initial data of small mass to  
scattering. For higher mass and sufficiently large $\alpha$, solitary 
waves seem to be observed for large times in accordance with the 
soliton resolution conjecture. However, for smaller values of 
$\alpha$, cusps can appear in finite time for such initial data. We also 
studied the formation of dispersive shock waves.

An interesting question raised by this study is to identify the 
parameter space $\alpha$, $\omega$ and $c$ for which smooth solitary 
waves exist. The fall-off behavior of these solutions should be 
proven. The orbital stability of these solitary waves is a 
question to be addressed also analytically. Of particular interest is 
the question of a blow-up already in the CH equation, and even more 
so in fCH, for which inital data a globally smooth solution in time 
can be expected, and for which data a blow-up in finite time is to be 
expected. The type of blow-up appears to be a gradient catastrophe, 
but this needs to be confirmed analytically. 

The formation of DSWs was shown numerically. In \cite{GK,AGK}, the 
onset of the oscillations as well as the boundary of the Whitham zone 
was conjectured to be asymptotically given by certain Painlev\'e 
transcendents. It is an interesting question whether there are 
fractional ODEs that play a role in this context for the fCH 
equation. It will be the subject of further research to address such 
questions.

\end{document}